\documentclass [12pt,oneside]{amsart}
\usepackage{amsmath}
\usepackage {amsfonts}
\usepackage {amsmath}
\usepackage {amsthm}
\usepackage {amssymb}
\usepackage {framed}
\usepackage {amsxtra}
\usepackage {enumerate}
\usepackage {graphicx}
\usepackage{color}
\usepackage{graphicx}
\usepackage[left=1 in, right=1 in, top=1.2 in, footskip=0.5 in, bottom=1in]{geometry}
\usepackage[font=small,skip=0pt]{caption}
\makeatletter

\theoremstyle{definition}
\newtheorem{df}{Definition} [section]

\theoremstyle{plain}
\newtheorem{thm}[df]{Theorem}

\newtheorem{conj}[df]{Conjecture}

\usepackage{tikz}
\usepackage{tikz-3dplot}
\usepackage{standalone}
\usetikzlibrary{matrix}
\usetikzlibrary{shapes,arrows,fit,calc,positioning,automata}

\title{On the Euclidean dimension of graphs}

\author{Jin Hyup Hong\\
Great Neck South High School, Great Neck, NY\\
and\\
Dan Ismailescu\\
Mathematics Department, Hofstra University, NY}

\begin{document}
\thispagestyle{empty}

\maketitle

\newpage
\thispagestyle{empty}

\begin{abstract}
The {\it Euclidean dimension a graph} $G$ is defined to be the smallest integer $d$ such that the vertices of $G$ can be located in $\mathbb{R}^d$ in such a way that two vertices are unit distance apart if and only if they are adjacent in $G$.
In this paper we determine the Euclidean dimension for twelve well known graphs. Five of these graphs, D\"{u}rer, Franklin, Desargues, Heawood and Tietze can be embedded in the plane, while the remaining graphs, Chv\'{a}tal, Goldner-Harrary, Herschel, Fritsch, Gr\"{o}tzsch, Hoffman and Soifer have Euclidean dimension $3$. We also present explicit embeddings for all these graphs.

\end{abstract}

\newpage
\pagenumbering{arabic}

\begin{section}{\bf History and previous work}

The {\it Euclidean dimension} of a graph $G=(V,E)$, denoted $dim(G)$ is the least integer $n$ such that
there exists a $1:1$ embedding $f:V\rightarrow \mathbb{R}^n$ for which $|f(u)-f(v)|=1$ if and only if $uv \in E$.

The concept was introduced by Erd\H{o}s, Harary and Tutte in their seminal paper \cite{EHT65}, where the authors determine
the Euclidean dimension for several classes of graphs.

For instance, they show that $dim(K_n) = n-1$, where $K_n$ is the complete graph on $n$ vertices.
Using a construction due to Lenz, they also compute the Euclidean dimension of $K_{m,n}$,
the complete bipartite graph with $m$ vertices in one class and $n$ vertices in the other.
\begin{thm}\cite{EHT65}
\vspace{-1cm}
\begin{align*}
&dim(K_{1,1})=1, \,\,dim(K_{1,n})=2\,\, \text{for}\,\, n\ge 2\\
&dim(K_{2,2})=2, \,\,dim(K_{2,n})=3\,\, \text{for}\,\, n\ge 3\\
&dim(K_{m,n})=4, \,\,\text{for all}\,\,\,\, n\ge m\ge  3.
\end{align*}
\end{thm}
Given two graphs $G_1=(V_1,E_1)$ and $G_2=(V_2,E_2)$ with disjoint vertex sets $V_1\cap V_2=\emptyset$, define
the {\it join} of these two graphs, denoted $G=G_1+G_2$, to be a graph $G=(V,E)$ such that $V=V_1\cup V_2$
and $E=E_1\cup E_2 \cup E'$ where $E'=\{v_1v_2 \,\,|\,\, v_1\in V_1 \,\,\text{and} \,\,v_2\in V_2\}$.
In other words, the join of two graphs is their graph union plus all possible edges joining each vertex of the first graph with
the vertices of the second graph.

For $n\ge 3$, let $W_{1,n}$ be the {\it wheel with $n$ spokes}, defined as $W_{1,n}=K_1+C_n$, the join of the one-vertex graph and the $n$ - cycle $C_n$. Erd\H{o}s, Harary and Tutte proved that $dim(W_{1,n})=3$ for all $n\neq 6$ and $dim(W_{1,6})=2$.

This is an interesting situation as it provides an instance where the Euclidean dimension of a graph is strictly smaller than
the Euclidean dimension of one of its subgraphs. Indeed, let $G=W_{1,6}$ which obviously has dimension $2$ as it can be seen from figure \ref{figwheel}. On the other hand, if one considers the subgraph $H=W_{1,6}-\{v_1v_2\}$, the wheel with a missing spoke, then $dim\,H\ge 3$ as one cannot embed $H$ in the plane and have $\|v_1-v_2\|\neq 1$.
\begin{figure}[!htb]
\centering
\includegraphics[scale=0.65]{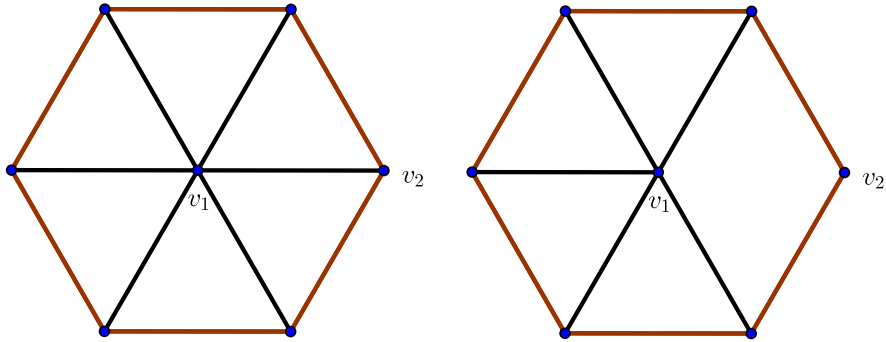}
\vspace{-1.cm}
\caption{An embedding of $W_{1,6}$ in the plane; $W_{1,6}\setminus \{v_1v_2\}$ cannot be embedded in $\mathbb{R}^2$.}
\label{figwheel}
\end{figure}
Buckley and Harary \cite{BH88} computed the Euclidean dimension for complete tripartite graphs $K_{m,n,p}$. Their result states that if $3\le m\le n\le p$ then $dim(K_{m,n,p})=6$; the cases when $\min(m,n,p)\le 2$ were settled as well.

 They also extended the results from \cite{EHT65} and determined the Euclidean dimension of the generalized wheel, $W_{m,n}=\overline{K_m}+C_n$, the join of the empty graph on $m$ vertices and the $n$ - cycle $C_n$.

Quite recently, Gervacio and Jos \cite{GJ08} determined the Euclidean dimension of the join of two cycles: their result states that
for all $m\ge 3$, $n\ge 3$ we have that $dim(C_m+C_n)=5$, except for $dim(C_4+C_4)=dim(C_5+C_5)=4$ and $dim(C_6+C_6)=6$.

In 2012, \v{Z}itnik, Horvat and Pisanski \cite{ZHP12} proved that all generalized Petersen graphs can be embedded as unit distance graphs in the plane. These graphs were introduced by Coxeter \cite{cox50} and studied again by Frucht, Graver and Watkins \cite{FGW71}.

In a series of recent preprints, Gerbracht \cite{G08,G09} studied the Euclidean dimension of symmetric trivalent graphs with up to 32 vertices and found $\mathbb{R}^2$ unit distance embeddings for many of them, several of which we mention here: M\"{o}bius-Kantor graph, the dodecahedral graph, Desargues' graph, the Nauru graph, the Levi (or Coxeter-Tutte) graph, the Dyck graph, Heawood graph.
We are going to revisit some of these graphs in the later sections.

It may appear that a lot is known about the Euclidean dimension of graphs. However, that is not the case.
Schaefer \cite{S13} proved that it is NP-hard to test whether the Euclidean dimension of a given graph is
at most a given value. The problem remains hard even for testing whether the Euclidean dimension is two.

The best general upper bound currently known is due to Maehara and R\"{odl}:
\begin{thm}\cite{MR90}\label{MR}
Let $G$ be a graph with maximum degree $\Delta$. Then $dim(G)\le 2\Delta$.
\end{thm}

This bound is in most cases very weak as it does not take into account the structure of the graph.
For instance, a complete bipartite graph can have arbitrarily large vertex degree but the Euclidean dimension
cannot be greater than 4.

To further illustrate how difficult the problem of finding the exact Euclidean dimension of a given
graph could be let us mention a famous example: the Heawood graph. This is the point-line incidence graph of the
finite projective plane of order two and it has $14$ vertices and $21$ edges. Chv\'{a}tal \cite{CKK72} suspected
that Heawood's graph cannot be embedded as a unit distance graph in the plane. Gerbracht \cite{G08} proved Chv\'{a}tal wrong
by finding eleven different planar unit distance embeddings of Heawood's graph. However, none of these constructions
is ``nice'' as the coordinates of the vertices depend on the roots of a polynomial of degree 79.
We are going to study Heawood's graph in one of the following sections.

In conclusion, there is no systematic method to determine the dimension of an arbitrary graph. This is on one hand unfortunate
but on the other hand it provides an intriguing list of open problems appropriate for a research project.
\end{section}

\begin{section}{\bf Our results}

In the sequel we determine the Euclidean dimension for twelve well-known graphs - see figure \ref{better}.
Five of these graphs have Euclidean dimension $2$: D\"{urer}, Franklin, Desargues, Heawood and Tietze. The remaining seven
graphs have Euclidean dimension $3$: Chv\'{a}tal, Goldner-Harary, Herschel, Fritsch, Gr\"{o}tzsch, Hoffman and Soifer.
For each of these graphs we provide explicit embeddings as well as a brief account on how these embedding were obtained.

Our initial interest in the problem was prompted by Gerbracht's paper \cite{G09} on Heawood's graph. We thought that such
a simple graph must have some more aesthetically pleasing embeddings than the ones Gerbracht found. We succeeded to find
an infinite family of axially symmetric embeddings of the Heawood graph. At that point we were not yet familiar with reference \cite{G08} since it was not publicly available. After we acquired a copy from Mr. Gerbracht himself, we noticed that he was interested exclusively in graphs of Euclidean dimension $2$. We thought it would be interesting to look in dimension $3$ as well.
\newpage
\begin{figure}[!htb]
\centering
\includegraphics[scale=0.33]{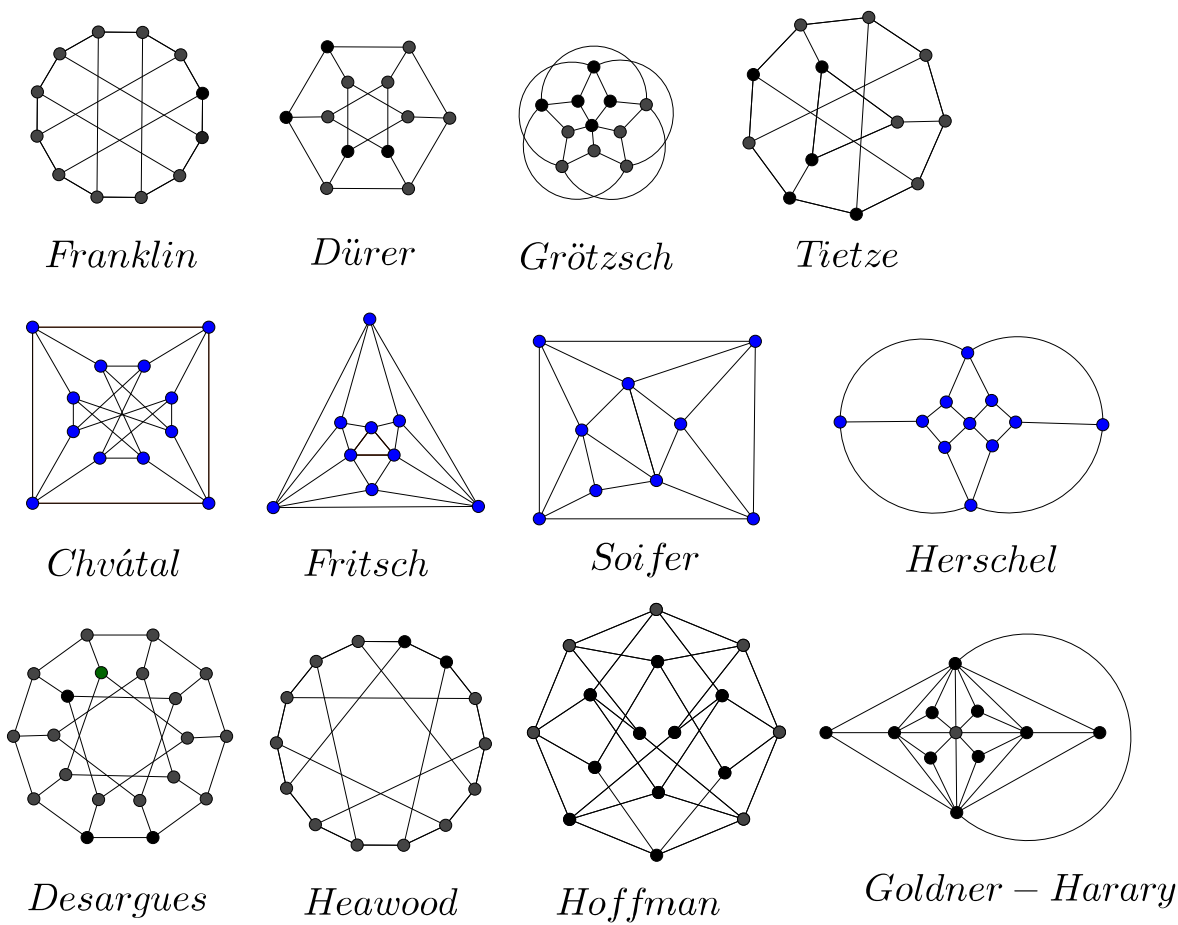}
\vspace{-2cm}
\caption{Twelve graphs}
\label{better}
\end{figure}
\vspace{-1cm}
Obviously, finding an embedding for a given graph basically requires solving a system of equations. The number of variables
equals the number of vertices of the graph times the dimension of the space where the embedding is attempted. The number of equations matches the number of edges of the particular graph. Solving this system by brute force is in most cases impossible.
We therefore use a few tricks to keep both the number of equations and the number of unknowns as small as possible.

Suppose we have two points in the plane $A_1(x_1,y_1)$ and $A_2(x_2,y_2)$ at unit distance from each other.
Thus, we have one equation $(x_1-x_2)^2+(y_1-y_2)^2=1$ and four unknowns, $x_1,y_1,x_2,y_2$. But we can very easily set $x_2=x_1+\cos{t}$ and $y_2=y_1+\sin{t}$ to get rid of the equation and decrease the number of variables from $4$ to $3$.
It may be argued that introducing trigonometric function makes the computation more difficult.

While this is true, we can circumvent this problem by using the substitutions $\cos{t}=(u^2-1)/(u^2+1)$ and $\sin{t}=2u/(u^2+1)$ where $u=\cot(t/2)$. It is true that in doing so we miss the value $t=0$, but in most cases that is not essential. We thus work with rational functions rather than trigonometric ones. We use the same idea in $\mathbb{R}^3$ as well.

Another useful approach is to exploit the rotational and/or axial symmetries of the graph. We use rotational symmetry to
produce embeddings for D\"{u}rer, Franklin, Desargues, Tietze and Gr\"{o}tzsch; we use the axial symmetry technique for Heawood.

Finally, we have to argue why each of the last seven graphs cannot be embedded in the plane, and have therefore Euclidean dimension $3$. The reason is what we call the {\it parallelogram impossibility condition} and we will describe it below.

Let us consider a very simple graph, the M\"{o}bius ladder on $6$ vertices - see figure \ref{mobius}.
\begin{figure}[!htb]
\centering
\includegraphics[scale=0.5]{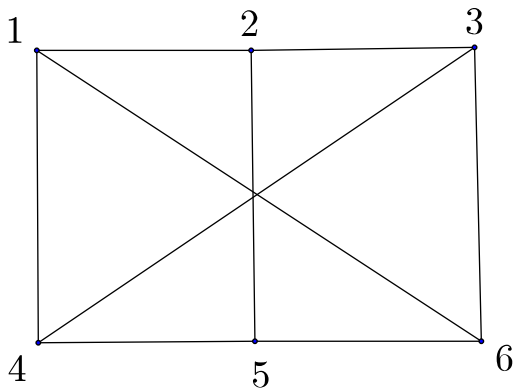}
\caption{M\"{o}bius Ladder Graph}
\label{mobius}
\end{figure}

We claim that this graph cannot be embedded in the plane as a unit distance graph. Suppose that there exists such an embedding.
Then quadrilateral $1254$ is a rhombus and therefore $\overrightarrow{14}=\overrightarrow{25}$, as vectors. The same argument holds for the rhombus $2365$; it follows that  $\overrightarrow{25}=\overrightarrow{36}$. But $1634$ must be a rhombus as well, since all sides have unit length. This gives that $\overrightarrow{14}=\overrightarrow{63}$. Combining the three equalities above,
it follows that  $\overrightarrow{63}=\overrightarrow{36}$, which means that vertices $3$ and $6$ coincide, contradiction.

Hence, the M\"{o}bius ladder has Euclidean dimension greater than $2$.
We use exactly the same argument for each of the seven graphs embedded in $\mathbb{R}^3$: we assume they can be embedded in the plane, list all the rhombi and then prove that two vertices must coincide.

\newpage
\hspace{-2cm}
\begin{tikzpicture}[scale=4.0,line width=1pt]
\tikzstyle{every node}=[draw=black,fill=yellow!50!white,thick,
  shape=circle,minimum size=2em,inner sep=0.4];
\draw (.5773502693, 0.)--(-.2886751347, .5000000000);
\draw (.5773502693, 0.)--(-.2886751347, -.5000000000);
\draw (.5773502693, 0.)--(.2886751347, .9574271080);
\draw (-.2886751347, .5000000000)--(-.2886751347, -.5000000000);
\draw (-.2886751347, .5000000000)--(-.9734937648, -.2287135539);
\draw (-.2886751347, -.5000000000)--(.6848186302, -.7287135539);
\draw (.2886751347, .5000000000)--(-.5773502693, 0.);
\draw (.2886751347, .5000000000)--(.2886751347, -.5000000000);
\draw (.2886751347, .5000000000)--(-.6848186302, .7287135539);
\draw (-.5773502693, 0.)--(.2886751347, -.5000000000);
\draw (-.5773502693, 0.)--(-.2886751347, -.9574271080);
\draw (.2886751347, -.5000000000)--(.9734937648, .2287135539);
\draw (.2886751347, .9574271080)--(-.6848186302, .7287135539);
\draw (.2886751347, .9574271080)--(.9734937648, .2287135539);
\draw (-.6848186302, .7287135539)--(-.9734937648, -.2287135539);
\draw (-.9734937648, -.2287135539)--(-.2886751347, -.9574271080);
\draw (-.2886751347, -.9574271080)--(.6848186302, -.7287135539);
\draw (.6848186302, -.7287135539)--(.9734937648, .2287135539);
\node at (.5773502693, 0) {1};
\node at (-.2886751347, .5000000000) {2};
\node at (-.2886751347, -.5000000000) {3};
\node at (.2886751347, .5000000000) {4};
\node at (-.5773502693, 0.) {5};
\node at (.2886751347, -.5000000000) {6};
\node at (.2886751347, .9574271080) {7};
\node at (-.6848186302, .7287135539) {8};
\node at (-.9734937648, -.2287135539) {9};
\node at (-.2886751347, -.9574271080) {10};
\node at (.6848186302, -.7287135539) {11};
\node at (.9734937648, .2287135539) {12};

\end{tikzpicture}
\tikzset{
    table/.style={
        matrix of nodes,
        row sep=-\pgflinewidth,
        column sep=-\pgflinewidth,
        nodes={
            rectangle,
            draw=black,
            align=center
        },
        minimum height=2.5em,
        text depth=1ex,
        text height=2ex,
        nodes in empty cells,
        every even row/.style={
            nodes={fill=gray!20}
        },
        column 1/.style={
            nodes={text width=2em,font=\bfseries}
        },
        row 1/.style={
            nodes={
                fill=black,
                text=white,
                font=\bfseries
            }
        }
    }
}
\begin{tikzpicture}

\matrix (first) [table,text width=10em]
{
& X & Y\\
1 & $\sqrt{3}/3$ & $0$ \\
2 & $-\sqrt{3}/6$ & $1/2$ \\
3 & $-\sqrt{3}/6$ & $-1/2$ \\
4 & $\sqrt{3}/6$ & $1/2$ \\
5 & $-\sqrt{3}/3$ & $0$ \\
6 & $\sqrt{3}/6$ & $-1/2$ \\
7 & $\sqrt{3}/6$ & $\sqrt{33}/6$ \\
8 & $\sqrt{3}/12-\sqrt{11}/4$ & $1/4+\sqrt{33}/12$ \\
9 & $-\sqrt{3}/12-\sqrt{11}/4$ & $1/4-\sqrt{33}/12$ \\
10 & $-\sqrt{3}/6$ & $-\sqrt{33}/6$ \\
11 & $-\sqrt{3}/12+\sqrt{11}/4$ & $-1/4-\sqrt{33}/12$ \\
12 & $\sqrt{3}/12+\sqrt{11}/4$ & $-1/4+\sqrt{33}/12$ \\
};
\end{tikzpicture}

\begin{figure}[!htb]
\centering
\caption{D\"{u}rer Graph}
\label{}
\end{figure}

In this case we use the $6$-fold symmetry of the graph.
In other words $A_4$ is obtained from $A_1$ after a $60^{\circ}$ rotation, $A_2$ is obtained from $A_4$ after a $60^{\circ}$ rotation, and so on. For short,
\begin{align*}
&A_1\longrightarrow A_4 \longrightarrow A_2 \longrightarrow A_5 \longrightarrow A_3 \longrightarrow A_6,\,\text{and}\\
&A_7\longrightarrow A_8 \longrightarrow A_9 \longrightarrow A_{10} \longrightarrow A_{11} \longrightarrow A_{12}.\\
\end{align*}

Since $|A_1-A_5|=1$ then we set $A_1(\sqrt{3}/3,0)$. The only thing left to do is to impose the conditions $|A_1-A_7|=1$ and $|A_7-A_{12}|=1$. These two equations give the coordinates of $A_7(\sqrt{3}/6,\sqrt{33}/6)$.


\newpage
\hspace{-2.2cm}
\begin{tikzpicture}[scale=4.3,line width=1pt]
\tikzstyle{every node}=[draw=black,fill=yellow!40!white,thick,
  shape=circle,minimum size=1.3em,inner sep=0.4];
    \draw (.9033103700, 0.)--(-0.0856297809, .1483151315);
    \draw (.9033103700, 0.)--(0.0372849659, .5000000000);
    \draw (0.372849659e-1, .5000000000)--(-.9516551849, .6483151315);
    \draw (.9033103700, 0.)--(-0.856297809e-1, -.1483151315);
    \draw (-.4516551849, .7822897280)--(-0.856297809e-1, -.1483151315);
    \draw (-0.856297809e-1, .1483151315)--(-.9516551849, .6483151315);
    \draw (-0.4516551849,-0.2177102725)--(-0.0856297809, -1.148315132);
    \draw (-0.856297809e-1, .1483151315)--(-.4516551849, -.7822897280);
    \draw (-.4516551849, -.7822897280)--(.1712595619, 0.);
    \draw (0.372849659e-1, .5000000000)--(1.037284966, .5000000000);
    \draw (.4143702191, -.2822897275)--(1.037284966, .5000000000);
    \draw (-.9516551849, .6483151315)--(-0.4516551849,-0.2177102725);
    \draw (-.4516551849, .7822897280)--(-0.4516551849,-0.2177102725);
    \draw (-.4516551849, .7822897280)--(.1712595619, 0.);
    \draw (-0.856297809e-1, -.1483151315)--(-0.856297809e-1, -1.148315132);
    \draw (-0.856297809e-1, -1.148315132)--(.4143702191, -.2822897275);
    \draw (-.4516551849, -.7822897280)--(.4143702191, -.2822897275);
    \draw (.1712595619, 0.)--(1.037284966, .5000000000);
\node at (.9033103700, 0.) {1};
\node at (-0.856297809e-1, .1483151315) {2};
\node at (0.372849659e-1, .5000000000) {3};
\node at (-.9516551849, .6483151315) {4};
\node at (-.4516551849, .7822897280) {5};
\node at (-0.856297809e-1, -.1483151315) {6};
\node at (-0.4516551849,-0.2177102725) {7};
\node at (-0.856297809e-1, -1.148315132) {8};
\node at (-.4516551849, -.7822897280) {9};
\node at (.1712595619, 0.) {10};
\node at (.4143702191, -.2822897275) {11};
\node at (1.037284966, .5000000000) {12};
\end{tikzpicture}
\begin{tikzpicture}
\tikzset{
    table/.style={
        matrix of nodes,
        row sep=-\pgflinewidth,
        column sep=-\pgflinewidth,
        nodes={
            rectangle,
            draw=black,
            align=center
        },
        minimum height=2.5em,
        text depth=1ex,
        text height=2ex,
        nodes in empty cells,
        every even row/.style={
            nodes={fill=gray!20}
        },
        column 1/.style={
            nodes={text width=2em,font=\bfseries}
        },
        row 1/.style={
            nodes={
                fill=black,
                text=white,
                font=\bfseries
            }
        }
    }
}
\matrix (first) [table,text width=11em]
{
& X & Y\\
1 & $3^{3/4}\sqrt{2}/6+\sqrt{3}/2-1/2$ & $0$ \\
2 & $\sqrt{3}/4-{1/4}-3^{3/4}\sqrt{2}/12$ & $-{3/4}+\sqrt{3}/4+3^{1/4}\sqrt{2}/4$ \\
3 & $3^{3/4}\sqrt{2}/6-1/2$ & $1/2$ \\
4 & $-1/4-\sqrt{3}/4-3^{3/4}\sqrt{2}/12$ & $-1/4+\sqrt{3}/4+3^{1/4}\sqrt{2}/4$ \\
5 & $-3^{3/4}\sqrt{2}/12-\sqrt{3}/4+1/4$ & $3^{1/4}\sqrt{2}/4+3/4-\sqrt{3}/4$ \\
6 & $\sqrt{3}/4-1/4-3^{3/4}\sqrt{2}/12$ & $3/4-\sqrt{3}/4-3^{1/4}\sqrt{2}/4$ \\
7 & $-3^{3/4}\sqrt{2}/12-\sqrt{3}/4+1/4$ & $3^{1/4}\sqrt{2}/4-\sqrt{3}/4-1/4$ \\
8 & $\sqrt{3}/4-1/4-3^{3/4}\sqrt{2}/12$ & $-3^{1/4}\sqrt{2}/4-\sqrt{3}/4-1/4$ \\
9 & $-3^{3/4}\sqrt{2}/12-\sqrt{3}/4+1/4$ & $-3^{1/4}\sqrt{2}/4-3/4+\sqrt{3}/4$ \\
10 & $-\sqrt{3}/2+1/2+3^{3/4}\sqrt{2}/6$ & $0$ \\
11 & $-3^{3/4}\sqrt{2}/12+\sqrt{3}/4+1/4$ & $-3^{1/4}\sqrt{2}/4+\sqrt{3}/4-1/4$ \\
12 & $1/2+3^{3/4}\sqrt{2}/6$ & $1/2$ \\
};
\end{tikzpicture}

\begin{figure}[!htb]
\centering
\caption{Franklin Graph}
\label{}
\end{figure}

In this case we use the $3$-fold symmetry of the graph. It is enough to define the first four vertices.
Set $A_1(a,0)$, $A_2=A_1+[\cos{\alpha},\sin{\alpha}]$, $A_3=A_2+[\cos{\beta},\sin{\beta}]$, $A_4=A_3+[\cos{\gamma},\sin{\gamma}]$,
then impose the conditions $|A_1-A_6|=1$ and $|A_3-A_{12}|=1$. We have two equations and four variables, hence there is plenty of freedom. We chose this solution since it looks quite nice.
\newpage
\hspace{-2.6cm}
\begin{tikzpicture}[scale=3.0,line width=1pt]
\tikzstyle{every node}=[draw=black,fill=yellow!50!white,thick,
  shape=circle,minimum size=1.3em,inner sep=0.4];
\draw (.6180339880, 0.)--(1.618033988, 0.);
\draw (.6180339880, 0.)--(-.1909830050, .5877852517);
\draw (.6180339880, 0.)--(-.1909830050, -.5877852517);
\draw (1.618033988, 0.)--(1.309016994, .9510565160);
\draw (1.618033988, 0.)--(1.309016994, -.9510565160);
\draw (.4999999993, .3632712636)--(1.309016994, .9510565160);
\draw (.4999999993, .3632712636)--(-.4999999993, .3632712636);
\draw (.4999999993, .3632712636)--(.1909830050, -.5877852517);
\draw (1.309016994, .9510565160)--(.4999999988, 1.538841768);
\draw (.1909830050, .5877852517)--(.4999999988, 1.538841768);
\draw (.1909830050, .5877852517)--(-.6180339880, 0.);
\draw (.1909830050, .5877852517)--(.4999999993, -.3632712636);
\draw (.4999999988, 1.538841768)--(-.4999999988, 1.538841768);
\draw (-.1909830050, .5877852517)--(-.4999999988, 1.538841768);
\draw (-.1909830050, .5877852517)--(-.4999999993, -.3632712636);
\draw (-.4999999988, 1.538841768)--(-1.309016994, .9510565160);
\draw (-.4999999993, .3632712636)--(-1.309016994, .9510565160);
\draw (-.4999999993, .3632712636)--(-.1909830050, -.5877852517);
\draw (-1.309016994, .9510565160)--(-1.618033988, 0.);
\draw (-.6180339880, 0.)--(-1.618033988, 0.);
\draw (-.6180339880, 0.)--(.1909830050, -.5877852517);
\draw (-1.618033988, 0.)--(-1.309016994, -.9510565160);
\draw (-.4999999993, -.3632712636)--(-1.309016994, -.9510565160);
\draw (-.4999999993, -.3632712636)--(.4999999993, -.3632712636);
\draw (-1.309016994, -.9510565160)--(-.4999999988, -1.538841768);
\draw (-.1909830050, -.5877852517)--(-.4999999988, -1.538841768);
\draw (-.4999999988, -1.538841768)--(.4999999988, -1.538841768);
\draw (.1909830050, -.5877852517)--(.4999999988, -1.538841768);
\draw (.4999999988, -1.538841768)--(1.309016994, -.9510565160);
\draw (.4999999993, -.3632712636)--(1.309016994, -.9510565160);

\node at (.6180339880, 0.) {1};
\node at (1.618033988, 0.) {2};
\node at (.4999999993, .3632712636) {3};
\node at (1.309016994, .9510565160) {4};
\node at (.1909830050, .5877852517) {5};
\node at (.4999999988, 1.538841768) {6};
\node at (-.1909830050, .5877852517) {7};
\node at (-.4999999988, 1.538841768) {8};
\node at (-.4999999993, .3632712636) {9};
\node at (-1.309016994, .9510565160) {10};
\node at (-.6180339880, 0.) {11};
\node at (-1.618033988, 0.) {12};
\node at (-.4999999993, -.3632712636) {13};
\node at (-1.309016994, -.9510565160) {14};
\node at (-.1909830050, -.5877852517) {15};
\node at (-.4999999988, -1.538841768) {16};
\node at (.1909830050, -.5877852517) {17};
\node at (.4999999988, -1.538841768) {18};
\node at (.4999999993, -.3632712636) {19};
\node at (1.309016994, -.9510565160) {20};

\end{tikzpicture}
\tikzset{
    table/.style={
        matrix of nodes,
        row sep=-\pgflinewidth,
        column sep=-\pgflinewidth,
        nodes={
            rectangle,
            draw=black,
            align=center
        },
        minimum height=2.0em,
        text depth=1ex,
        text height=1.3ex,
        nodes in empty cells,
        every even row/.style={
            nodes={fill=gray!20}
        },
        column 1/.style={
            nodes={text width=2em,font=\bfseries}
        },
        row 1/.style={
            nodes={
                fill=black,
                text=white,
                font=\bfseries
            }
        }
    }
}
\begin{tikzpicture}

\matrix (first) [table,text width=9em]
{
& X & Y\\
1 & $\sqrt{5}/2-1/2$ & $0$ \\
2 & $1/2+\sqrt{5}/2$ & $0$ \\
3 & {\footnotesize$(\sqrt{5}/2-1/2)\cos(\pi/5)$} & {\footnotesize$(\sqrt{5}/2-1/2)\sin(\pi/5)$} \\
4 & {\footnotesize$(\sqrt{5}/2+1/2)\cos(\pi/5)$} & {\footnotesize$(\sqrt{5}/2+1/2)\sin(\pi/5)$} \\
5 & {\footnotesize $(\sqrt{5}/2-1/2)\cos(2\pi/5)$} & {\footnotesize$(\sqrt{5}/2-1/2)\sin(2\pi/5)$} \\
6 & {\footnotesize $(\sqrt{5}/2+1/2)\cos(2\pi/5)$} & {\footnotesize$(\sqrt{5}/2+1/2)\sin(2\pi/5)$} \\
7 & {\scriptsize $(-\sqrt{5}/2+1/2)\cos(2\pi/5)$} & {\footnotesize$(\sqrt{5}/2-1/2)\sin(2\pi/5)$} \\
8 & {\scriptsize $-(\sqrt{5}/2+1/2)\cos(2\pi/5)$} & {\footnotesize$(\sqrt{5}/2+1/2)\sin(2\pi/5)$} \\
9 & {\scriptsize $-(\sqrt{5}/2-1/2)\cos(\pi/5)$} & {\footnotesize $(\sqrt{5}/2-1/2)\sin(\pi/5)$} \\
10 & {\scriptsize $-(1/2+\sqrt{5}/2)\cos(\pi/5)$} & {\footnotesize $(1/2+\sqrt{5}/2)\sin(\pi/5)$} \\
11 & $-\sqrt{5}/2+1/2$ & $0$ \\
12 & $-1/2-\sqrt{5}/2$ & $0$ \\
13 & {\scriptsize $-(\sqrt{5}/2-1/2)\cos(\pi/5)$} & {\scriptsize $-(\sqrt{5}/2-1/2)\sin(\pi/5)$} \\
14 & {\scriptsize $-(1/2+\sqrt{5}/2)\cos(\pi/5)$} & {\scriptsize $-(1/2+\sqrt{5}/2)\sin(\pi/5)$} \\
15 & {\scriptsize $-(\sqrt{5}/2-1/2)\cos(2\pi/5)$} & {\scriptsize $-(\sqrt{5}/2-1/2)\sin(2\pi/5)$} \\
16 & {\scriptsize $-(1/2+\sqrt{5}/2)\cos(2\pi/5)$} & {\scriptsize $-(1/2+\sqrt{5}/2)\sin(2\pi/5)$} \\
17 & {\footnotesize $(\sqrt{5}/2-1/2)\cos((2\pi/5)$} & {\scriptsize $-(\sqrt{5}/2-1/2)\sin(2\pi/5)$} \\
18 & {\footnotesize $(1/2+\sqrt{5}/2)\cos(2\pi/5)$} & {\scriptsize $-(1/2+\sqrt{5}/2)\sin(2\pi/5)$} \\
19 & {\scriptsize $(\sqrt{5}/2-1/2)\cos(\pi/5)$} & {\scriptsize $-(\sqrt{5}/2-1/2)\sin(\pi/5)$} \\
20 & {\footnotesize $(1/2+\sqrt{5}/2)\cos(\pi/5)$} & {\scriptsize $-(1/2+\sqrt{5}/2)\sin(\pi/5)$} \\
};
\end{tikzpicture}

\begin{figure}[!htb]
\centering
\caption{Desargues Graph}
\label{}
\end{figure}
We use the $10$-fold rotational symmetry. It is sufficient to locate two vertices, since the rest can be obtained
via successive $36^{\circ}$ rotations. The same embedding appears in \cite{G08} but be rediscovered it independently.

\newpage
\hspace{-2.0cm}
\begin{tikzpicture}[scale=4,line width=1pt]
\node at (-0.55,0) {1};
\node at (0.55,0) {2};
\node at (1.15416457148644122770821426387, .804343497512308899008221633856) {3};
\node at (.104440118855331767667374281465, .601355708299205927066976097201) {4};
\node at (-.300414997537589923253483411073, 1.50455980435754223282095036117) {5};
\node at (-1.14416457148644122770821426387, .804343497512308899008221633856) {6};
\node at (-.104440118855331767667374281465, .601355708299205927066976097201) {7};
\node at (.300414997537589923253483411073, 1.50455980435754223282095036117) {8};
\node at (-0.6, 1.95959179422654247855782725978) {9};
\node at (0.6, 1.95959179422654247855782725978) {\footnotesize $10$};
\node at (0.4, 0.979795897113271239278913629884) {\footnotesize $11$};
\node at (-0.6, 1.57979589711327123927891362989) {\footnotesize $12$};
\node at (0.6, 1.57979589711327123927891362989) {\footnotesize $13$};
\node at (-0.4, 0.979795897113271239278913629884) {\footnotesize $14$};
\tikzstyle{every node}=[draw=red,fill=red!20!,shape=circle,scale=0.25]
\draw (-1/2,0)--(1/2,0);
\draw (-1/2,0)--(-1.09416457148644122770821426387, .804343497512308899008221633856);
\draw (-1/2,0)--(-3/10, 0.979795897113271239278913629884);
\draw (1/2,0)--(1.09416457148644122770821426387, .804343497512308899008221633856);
\draw (1/2,0)--(3/10, 0.979795897113271239278913629884);
\draw (1.09416457148644122770821426387, .804343497512308899008221633856)--(.104440118855331767667374281465, .661355708299205927066976097201);
\draw (1.09416457148644122770821426387, .804343497512308899008221633856)--(.400414997537589923253483411073, 1.52455980435754223282095036117);
\draw (.104440118855331767667374281465, .661355708299205927066976097201)--(-.400414997537589923253483411073, 1.52455980435754223282095036117);
\draw (.104440118855331767667374281465, .661355708299205927066976097201)--(1/2, 1.57979589711327123927891362989);
\draw (-.400414997537589923253483411073, 1.52455980435754223282095036117)--(-1.09416457148644122770821426387, .804343497512308899008221633856);
\draw (-.400414997537589923253483411073, 1.52455980435754223282095036117)--(1/2, 1.95959179422654247855782725978);
\draw (-1.09416457148644122770821426387, .804343497512308899008221633856)--(-.104440118855331767667374281465, .661355708299205927066976097201);
\draw (-.104440118855331767667374281465, .661355708299205927066976097201)--(.400414997537589923253483411073, 1.52455980435754223282095036117);
\draw (-.104440118855331767667374281465, .661355708299205927066976097201)--(-1/2, 1.57979589711327123927891362989);
\draw (.400414997537589923253483411073, 1.52455980435754223282095036117)--(-1/2, 1.95959179422654247855782725978);
\draw (-1/2, 1.95959179422654247855782725978)--(1/2, 1.95959179422654247855782725978);
\draw (-1/2, 1.95959179422654247855782725978)--(-3/10, 0.979795897113271239278913629884);
\draw (1/2, 1.95959179422654247855782725978)--(3/10, 0.979795897113271239278913629884);
\draw (3/10, 0.979795897113271239278913629884)--(-1/2, 1.57979589711327123927891362989);
\draw (-1/2, 1.57979589711327123927891362989)--(1/2, 1.57979589711327123927891362989);
\draw (1/2, 1.57979589711327123927891362989)--(-3/10, 0.979795897113271239278913629884);
\node at (-1/2,0) {};
\node at (1/2,0) {};
\node at (1.09416457148644122770821426387, .804343497512308899008221633856) {};
\node at (.104440118855331767667374281465, .661355708299205927066976097201) {};
\node at (-.400414997537589923253483411073, 1.52455980435754223282095036117) {};
\node at (-1.09416457148644122770821426387, .804343497512308899008221633856) {};
\node at (-.104440118855331767667374281465, .661355708299205927066976097201) {};
\node at (.400414997537589923253483411073, 1.52455980435754223282095036117) {};
\node at (-1/2, 1.95959179422654247855782725978) {};
\node at (1/2, 1.95959179422654247855782725978) {};
\node at (3/10, 0.979795897113271239278913629884) {};
\node at (-1/2, 1.57979589711327123927891362989) {};
\node at (1/2, 1.57979589711327123927891362989) {};
\node at (-3/10, 0.979795897113271239278913629884) {};
\end{tikzpicture}
\begin{tikzpicture}
\tikzset{
    table/.style={
        matrix of nodes,
        row sep=-\pgflinewidth,
        column sep=-\pgflinewidth,
        nodes={
            rectangle,
            draw=black,
            align=center
        },
        minimum height=2.5em,
        text depth=1ex,
        text height=2ex,
        nodes in empty cells,
        every even row/.style={
            nodes={fill=gray!20}
        },
        column 1/.style={
            nodes={text width=2em,font=\bfseries}
        },
        row 1/.style={
            nodes={
                fill=black,
                text=white,
                font=\bfseries
            }
        }
    }
}

\matrix (first) [table,text width=7em]
{
& X & Y \\
1 & $-1/2$ & $0$ \\
2 & $1/2$ & $0$ \\
3 & $1.094164572$ & $0.8043434978$ \\
4 & $0.1044401192$ & $0.6613557099$ \\
5 & $-0.4004149975$ & $1.524559804$ \\
6 & $-1.094164572$ & $0.8043434978$ \\
7 & $-0.1044401192$ & $0.6613557099$ \\
8 & $0.4004149975$ & $1.524559804$ \\
9 & $-1/2$ & $1.959591792$ \\
10 & $1/2$ & $1.959591792$ \\
11 & $3/10$ & $2\sqrt{6}/5$ \\
12 & $1/2$ & $1.579795904$ \\
13 & $1/2$ & $1.579795904$ \\
14 & $-3/10$ & $2\sqrt{6}/5$ \\
};
\end{tikzpicture}
\begin{figure}[!htb]
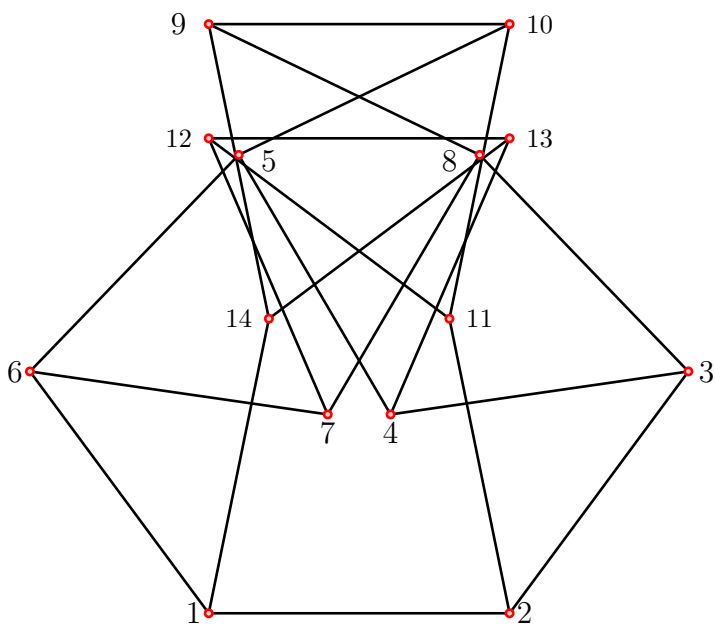

\centering
\caption{Heawood Graph}
\label{}
\end{figure}

The embbeding is axially symmetric. This was a tough one but we got a much better looking embedding
than Gerbracht. Some coordinates are given in numerical form since they depend of the roots of a polynomial of degree $14$
with coefficients in $\mathbb{Z}(\sqrt{6})$
{\tiny{
\begin{align*}
&216000x^{14}+(-847728\sqrt{6}-1052352)x^{13}+(4792176\sqrt{6}+9023184)x^{12}+(-13587312\sqrt{6}-51517008)x^{11}+(49227504\sqrt{6}+112815936)x^{10}+\\
&+(-117320465\sqrt{6}-171064560)x^9+(67247905\sqrt{6}+500137520)x^8+(-231804550\sqrt{6}-165859200)x^7+(96507350\sqrt{6}+411138400)x^6+\\
&(-94161060\sqrt{6}-273683040)x^5+(63781020\sqrt{6}+228901680)x^4+(-118255138\sqrt{6}+13784208)x^3+(-25457854\sqrt{6}+215805664)x^2+\\
&+(-31967747\sqrt{6}+25007952)x-1194101\sqrt{6}+11113616.
\end{align*}}}
\hspace{-2.5cm}
\begin{tikzpicture}[scale=4.7,line width=1pt]
\tikzstyle{every node}=[draw=black,fill=yellow!40!white,thick,
  shape=circle,minimum size=1.3em,inner sep=0.4];

\draw (.5773502693, 0.)--(-.4226497307, 0.);
\draw (.5773502693, 0.)--(-.2886751347, .5000000000);
\draw (.5773502693, 0.)--(-.2886751347, -.5000000000);
\draw (-.4226497307, 0.)--(.3596399968, .6229147468);
\draw (-.4226497307, 0.)--(-.5709648622, -.9889401510);
\draw (-.4226497307, 0.)--(.3596399968, -.6229147468);
\draw (-.4226497307, 0.)--(-.5709648622, .9889401510);
\draw (.3596399968, .6229147468)--(.2113248653, -.3660254040);
\draw (.3596399968, .6229147468)--(1.141929724, 0.);
\draw (.3596399968, .6229147468)--(-.5709648622, .9889401510);
\draw (-.5709648622, -.9889401510)--(.2113248653, -.3660254040);
\draw (-.5709648622, -.9889401510)--(-.7192799937, 0.);
\draw (-.5709648622, -.9889401510)--(.3596399968, -.6229147468);
\draw (-.2886751347, .5000000000)--(.2113248653, -.3660254040);
\draw (-.2886751347, .5000000000)--(-.2886751347, -.5000000000);
\draw (.2113248653, -.3660254040)--(-.7192799937, 0.);
\draw (.2113248653, -.3660254040)--(1.141929724, 0.);
\draw (-.7192799937, 0.)--(.2113248653, .3660254040);
\draw (-.7192799937, 0.)--(-.5709648622, .9889401510);
\draw (1.141929724, 0.)--(.2113248653, .3660254040);
\draw (1.141929724, 0.)--(.3596399968, -.6229147468);
\draw (-.2886751347, -.5000000000)--(.2113248653, .3660254040);
\draw (.2113248653, .3660254040)--(.3596399968, -.6229147468);
\draw (.2113248653, .3660254040)--(-.5709648622, .9889401510);
\node at (.5773502693, 0.) {1};
\node at (-.4226497307, 0.) {2};
\node at (.3596399968, .6229147468) {3};
\node at (-.5709648622, -.9889401510) {4};
\node at (-.2886751347, .5000000000) {5};
\node at (.2113248653, -.3660254040) {6};
\node at (-.7192799937, 0.) {7};
\node at (1.141929724, 0.) {8};
\node at (-.2886751347, -.5000000000) {9};
\node at (.2113248653, .3660254040) {10};
\node at (.3596399968, -.6229147468) {11};
\node at (-.5709648622, .9889401510) {12};
\end{tikzpicture}
\begin{tikzpicture}
\tikzset{
    table/.style={
        matrix of nodes,
        row sep=-\pgflinewidth,
        column sep=-\pgflinewidth,
        nodes={
            rectangle,
            draw=black,
            align=center
        },
        minimum height=2.5em,
        text depth=1ex,
        text height=2ex,
        nodes in empty cells,
        every even row/.style={
            nodes={fill=gray!20}
        },
        column 1/.style={
            nodes={text width=2em,font=\bfseries}
        },
        row 1/.style={
            nodes={
                fill=black,
                text=white,
                font=\bfseries
            }
        }
    }
}

\matrix (first) [table,text width=10em]
{
& X & Y\\
1 & $\sqrt{3}/3$ & $0$ \\
2 & $\sqrt{3}/3-1$ & $0$ \\
3 & $\sqrt{3}/12+3^{1/4}\sqrt{2}/4-1/4$ & $3^{3/4}\sqrt{2}/4-\sqrt{3}/4+1/4$ \\
4 & $\sqrt{3}/12-3^{1/4}\sqrt{2}/4-1/4$ & $1/4-3^{3/4}\sqrt{2}-\sqrt{3}/4$ \\
5 & $-\sqrt{3}/6$ & $1/2$ \\
6 & $1/2-\sqrt{3}/6$ & $1/2-\sqrt{3}$ \\
7 & $1/2-\sqrt{3}/6-3^{1/4}\sqrt{2}/2$ & $0$ \\
8 & $1/2+3^{1/4}\sqrt{2}/2-\sqrt{3}/6$ & $0$ \\
9 & $-\sqrt{3}/6$ & $-1/2$ \\
10 & $1/2-\sqrt{3}/6$ & $\sqrt{3}/2-1/2$ \\
11 & $3^{1/4}\sqrt{2}/4+\sqrt{3}/12-1/4$ & $\sqrt{3}/4-3^{3/4}\sqrt{2}/4-1/4$ \\
12 & $\sqrt{3}/12-3^{1/4}\sqrt{2}/4-1/4$ & $3^{3/4}\sqrt{2}/4+\sqrt{3}/4-1/4$ \\
};
\end{tikzpicture}
\begin{figure}[!htb]
\centering
\caption{Tietze Graph}
\label{}
\end{figure}

For constructing this graph we exploit its natural $3$-fold rotational symmetry.
It is sufficient to locate the first four vertices since the rest can be then obtained via
$120^{\circ}$ rotations.
Since $|A_1-A_5|=1$, we chose $A_1(\sqrt{3}/3,0)$. Since $|A_1-A_2|=1$, we select
$A_2= A_1+[(t_2^2-1)/(t_2^2+1), 2t_2/(t_2^2+1)]$. Since $|A_2-A_3|=1$ and $|A_2-A_4|=1$ we let
$A_3 = A_2+[(t_3^2-1)/(t_3^2+1), 2t_3/(t_3^2+1)]$ and $A_4= A_2+[(t_4^2-1)/(t_4^2+1), 2t_4/(t_4^2+1)]$.
We impose the conditions $|A_3-A_8|=1$ and $|A_3-A_{12}|=1$ to account of the remaining edges.
We have a system of two equations and three unknowns that happens to have nice solutions.
The one we chose is $t_2 = 0,\, t_3 = 1+3^{1/4}\,\sqrt{2},\, t_4 = 1-3^{1/4}\,\sqrt{2}$.
This choice generates the graph presented in the figure.
\newpage

\vspace{-5cm}
\hspace{-3.0cm}
\tdplotsetmaincoords{-75}{146}
\begin{tikzpicture}[tdplot_main_coords,
            cube/.style={very thick,black},
			grid/.style={gray,very thin,opacity=0.05},
			axis/.style={->,blue,thick},
            line width=0.7pt,scale=5.0]
\draw[axis] (0,0,0) -- (1.0,0,0) node[anchor=north]{$x$};
	\draw[axis] (0,0,0) -- (0,1.2,0) node[anchor=west]{$y$};
	\draw[axis] (0,0,0) -- (0,0,1.2) node[anchor=west]{$z$};
\node at (.644275050950098216971912173198, .291392619541484520778353627239, 1.045661171324115547230637687478) {1};
\node at (.550000000000000000000000000000, 0., 0.05) {2};
\node at (-.353553390593273762200422181052, .253553390593273762200422181052, 1.00879532511286756128183189395) {3};
\node at (-.291392619541484520778353627239, .644275050950098216971912173198, 1.045661171324115547230637687478) {4};
\node at (0.03, .400000000000000000000000000000, 0.05) {5};
\node at (-.330553390593273762200422181052, -.320553390593273762200422181052, 0.803879532511286756128183189395) {6};
\node at (-.644275050950098216971912173198, -.291392619541484520778353627239, 1.005661171324115547230637687478) {7};
\node at (-.500000000000000000000000000000, 0., -0.05) {8};
\node at (.353553390593273762200422181052, -.353553390593273762200422181052, 1.123879532511286756128183189395) {9};
\node at (.291392619541484520778353627239, -.644275050950098216971912173198, 1.045661171324115547230637687478) {10};
\node at (0., -.500000000000000000000000000000, -0.05) {11};
\node at (.403553390593273762200422181052, .353553390593273762200422181052, 1.023879532511286756128183189395) {12};
\draw (.644275050950098216971912173198, .291392619541484520778353627239, .945661171324115547230637687478)--(.5,0,0);
\draw (.644275050950098216971912173198, .291392619541484520778353627239, .945661171324115547230637687478)--(-.353553390593273762200422181052, .353553390593273762200422181052, .923879532511286756128183189395);
\draw (.644275050950098216971912173198, .291392619541484520778353627239, .945661171324115547230637687478)--(-.291392619541484520778353627239, .644275050950098216971912173198, .945661171324115547230637687478);
\draw (.644275050950098216971912173198, .291392619541484520778353627239, .945661171324115547230637687478)--(.291392619541484520778353627239, -.644275050950098216971912173198, .945661171324115547230637687478);
\draw (.5,0,0)--(-.500000000000000000000000000000, 0., 0);
\draw (.5,0,0)--(.353553390593273762200422181052, -.353553390593273762200422181052, .923879532511286756128183189395);
\draw (.5,0,0)--(.353553390593273762200422181052, .353553390593273762200422181052, .923879532511286756128183189395);
\draw (-.353553390593273762200422181052, .353553390593273762200422181052, .923879532511286756128183189395)--(0,0.5,0);
\draw (-.353553390593273762200422181052, .353553390593273762200422181052, .923879532511286756128183189395)--(-.500000000000000000000000000000, 0., 0);
\draw (-.353553390593273762200422181052, .353553390593273762200422181052, .923879532511286756128183189395)--(.353553390593273762200422181052, -.353553390593273762200422181052, .923879532511286756128183189395);
\draw (-.291392619541484520778353627239, .644275050950098216971912173198, .945661171324115547230637687478)--(0., .500000000000000000000000000000, 0);
\draw (-.291392619541484520778353627239, .644275050950098216971912173198, .945661171324115547230637687478)--(-.353553390593273762200422181052, -.353553390593273762200422181052, .923879532511286756128183189395);
\draw (-.291392619541484520778353627239, .644275050950098216971912173198, .945661171324115547230637687478)--(-.644275050950098216971912173198, -.291392619541484520778353627239, .945661171324115547230637687478);
\draw (0., .500000000000000000000000000000, 0)--(0., -.500000000000000000000000000000, 0);
\draw (0., .500000000000000000000000000000, 0)--(.353553390593273762200422181052, .353553390593273762200422181052, .923879532511286756128183189395);
\draw (-.353553390593273762200422181052, -.353553390593273762200422181052, .923879532511286756128183189395)--(-.500000000000000000000000000000, 0., 0);
\draw (-.353553390593273762200422181052, -.353553390593273762200422181052, .923879532511286756128183189395)--(0., -.500000000000000000000000000000, 0);
\draw (-.353553390593273762200422181052, -.353553390593273762200422181052, .923879532511286756128183189395)--(.353553390593273762200422181052, .353553390593273762200422181052, .923879532511286756128183189395);
\draw (-.644275050950098216971912173198, -.291392619541484520778353627239, .945661171324115547230637687478)--(-.500000000000000000000000000000, 0., 0);
\draw (-.644275050950098216971912173198, -.291392619541484520778353627239, .945661171324115547230637687478)--(.353553390593273762200422181052, -.353553390593273762200422181052, .923879532511286756128183189395);
\draw (-.644275050950098216971912173198, -.291392619541484520778353627239, .945661171324115547230637687478)--(.291392619541484520778353627239, -.644275050950098216971912173198, .945661171324115547230637687478);
\draw (.353553390593273762200422181052, -.353553390593273762200422181052, .923879532511286756128183189395)--(0., -.500000000000000000000000000000, 0);
\draw (.291392619541484520778353627239, -.644275050950098216971912173198, .945661171324115547230637687478)--(0., -.500000000000000000000000000000, 0);
\draw (.291392619541484520778353627239, -.644275050950098216971912173198, .945661171324115547230637687478)--(.353553390593273762200422181052, .353553390593273762200422181052, .923879532511286756128183189395);
\tikzstyle{every node}=[draw=red,fill=red!20!,shape=circle,scale=0.25]
\node at (.644275050950098216971912173198, .291392619541484520778353627239, .945661171324115547230637687478) {};
\node at (.500000000000000000000000000000, 0., 0) {};
\node at (-.353553390593273762200422181052, .353553390593273762200422181052, .923879532511286756128183189395) {};
\node at (-.291392619541484520778353627239, .644275050950098216971912173198, .945661171324115547230637687478) {};
\node at (0., .500000000000000000000000000000, 0) {};
\node at (-.353553390593273762200422181052, -.353553390593273762200422181052, .923879532511286756128183189395) {};
\node at (-.644275050950098216971912173198, -.291392619541484520778353627239, .945661171324115547230637687478) {};
\node at (-.500000000000000000000000000000, 0., 0) {};
\node at (.353553390593273762200422181052, -.353553390593273762200422181052, .923879532511286756128183189395) {};
\node at (.291392619541484520778353627239, -.644275050950098216971912173198, .945661171324115547230637687478) {};
\node at (0., -.500000000000000000000000000000, 0) {};
\node at (.353553390593273762200422181052, .353553390593273762200422181052, .923879532511286756128183189395) {};
\end{tikzpicture}
\begin{tikzpicture}[scale=0.5]
\tikzset{
    table/.style={
        matrix of nodes,
        row sep=-\pgflinewidth,
        column sep=-\pgflinewidth,
        nodes={
            rectangle,
            draw=black,
            align=center
        },
        minimum height=1.5em,
        text depth=1ex,
        text height=2ex,
        nodes in empty cells,
        every even row/.style={
            nodes={fill=gray!20}
        },
        column 1/.style={
            nodes={text width=1em,font=\bfseries}
        },
        row 1/.style={
            nodes={
                fill=black,
                text=white,
                font=\bfseries
            }
        }
    }
}
 \pgfgettransform\mytrafo
\matrix (first) [table,text width=8.0em,execute at begin cell=\pgfsettransform\mytrafo]
{
& X & Y & Z\\
1 & $(b^{2}-1)/\sqrt{2}(b^{2}+1)$ & $\sqrt{2}b/(b^{2}+1)$ & $c$ \\
2 & $1/2$ & $0$ & $0$\\
3 & $-\sqrt{2}/4$ & $\sqrt{2}/4$ & $a/2+a/2\sqrt{2}$ \\
4 & $-b\sqrt{2}/(b^{2}+1)$ & ${b^{2}-1}/\sqrt{2}(b^{2}+1)$ & $c$\\
5 & $0$ & $1/2$ & $0$ \\
6 & $-\sqrt{2}/4$ & $-\sqrt{2}/4$ & $a/2+a/2\sqrt{2}$ \\
7 & {\small $-(b^{2}-1)/\sqrt{2}(b^{2}+1)$} & $-\sqrt{2}b/(b^{2}+1)$ & $c$ \\
8 & $-1/2$ & $0$ & $0$ \\
9 & $\sqrt{2}/4$ & $\sqrt{2}/4$ & $a/2+a/2\sqrt{2}$ \\
10 & $b\sqrt{2}/(b^{2}+1)$ & {\small $-(b^{2}-1)/\sqrt{2}(b^{2}+1)$} & $c$\\
11 & $0$ & $-1/2$ & $0$ \\
12 & $\sqrt{2}/4$ & $\sqrt{2}/4$ & $a/2+a/2\sqrt{2}$ \\
};
\end{tikzpicture}

\begin{figure}[!htb]
\centering
\caption{Chv\'atal Graph: $a$=$\sqrt{4-2\sqrt{2}}$, $b$=$\sqrt{(2\sqrt{2}-1+4c^{2})/(\sqrt{2}+1-4c^{2})}$, $c$=$\sqrt{2}\sqrt{2-\sqrt{2}}/4+
\sqrt{2-\sqrt{2}}/4+\sqrt{\sqrt{2}+4\sqrt{2}\sqrt{2-\sqrt{2}}-2}/4$}
\label{}
\end{figure}

The idea is to use the $4$-fold rotational symmetry of the graph. We want $A_1A_4A_7A_{10}$ to be a unit square, and $A_2A_5A_8A_{11}$ and $A_3A_6A_9A_12$ are squares with diagonal equal to $1$. These three squares lie in planes parallel
to the $xy$ plane. The conditions $|A_1-A_2|=1$ and $|A_1-A_3|=1$ give the above solution.
\newpage
\hspace{-3.0cm}
\tdplotsetmaincoords{67}{163}
\begin{tikzpicture}[tdplot_main_coords,
            cube/.style={very thick,black},
			grid/.style={gray,very thin,opacity=0.05},
			axis/.style={->,blue,very thick},
            line width=1.0pt,scale=3.6]
\foreach \x in {-0.8,-0.4,...,1.2}
\foreach \y in {-0.8,-0.4,...,1.2}
\foreach \z in {-0.8,-0.4,...,1.2}
	{
			\draw[grid] (\x,-0.8,-0.8) -- (\x,1.2,-0.8);
			\draw[grid] (-0.8,\y,-0.8) -- (1.2,\y,-0.8);
            \draw[grid] (-0.8,-0.8,\z) -- (1.2,-0.8,\z);
            \draw[grid] (-0.8,-0.8,\z) -- (-0.8,1.2,\z);
           \draw[grid] (\x,-0.8,-0.8) -- (\x,-0.8,1.2);
            \draw[grid] (-0.8,\y,-0.8) -- (-0.8,\y,1.2);
		}
\draw[axis] (0,0,0) -- (1.6,0,0) node[anchor=south]{$x$};
	\draw[axis] (0,0,0) -- (0,2.8,0) node[anchor=west]{$y$};
	\draw[axis] (0,0,0) -- (0,0,1.6) node[anchor=west]{$z$};
\draw (0., 0., 0.)--(1., 0., 0.);
\draw (0., 0., 0.)--(.5000000000, -.6735753142, -.5443310539);
\draw (0., 0., 0.)--(.5000000000, .2886751347, .8164965809);
\draw (0., 0., 0.)--(.5000000000, .2886751347, -.8164965809);
\draw (0., 0., 0.)--(.5000000000, .8660254040, 0.);
\draw (0., 0., 0.)--(-.3888888889, -.1603750748, -.9072184234);
\draw (0., 0., 0.)--(-.3888888889, .8018753741, .4536092117);
\draw (0., 0., 0.)--(-.3333333333, .7698003590, -.5443310539);
\draw (1., 0., 0.)--(.5000000000, -.6735753142, -.5443310539);
\draw (1., 0., 0.)--(.5000000000, .2886751347, .8164965809);
\draw (1., 0., 0.)--(.5000000000, .2886751347, -.8164965809);
\draw (1., 0., 0.)--(.5000000000, .8660254040, 0.);
\draw (1., 0., 0.)--(1.333333333, .7698003590, -.5443310539);
\draw (.5000000000, -.6735753142, -.5443310539)--(.5000000000, .2886751347, -.8164965809);
\draw (.5000000000, .2886751347, .8164965809)--(.5000000000, .8660254040, 0.);
\draw (.5000000000, .2886751347, -.8164965809)--(.5000000000, .8660254040, 0.);
\draw (.5000000000, .2886751347, -.8164965809)--(1.333333333, .7698003590, -.5443310539);
\draw (.5000000000, .2886751347, -.8164965809)--(-.3888888889, -.1603750748, -.9072184234);
\draw (.5000000000, .2886751347, -.8164965809)--(-.3333333333, .7698003590, -.5443310539);
\draw (.5000000000, .2886751347, -.8164965809)--(.4444444444, 1.283000598, -.9072184234);
\draw (.5000000000, .8660254040, 0.)--(1.333333333, .7698003590, -.5443310539);
\draw (.5000000000, .8660254040, 0.)--(-.3888888889, .8018753741, .4536092117);
\draw (.5000000000, .8660254040, 0.)--(-.3333333333, .7698003590, -.5443310539);
\draw (.5000000000, .8660254040, 0.)--(.4444444444, 1.283000598, -.9072184234);
\draw (-.3888888889, -.1603750748, -.9072184234)--(-.3333333333, .7698003590, -.5443310539);
\draw (-.3888888889, .8018753741, .4536092117)--(-.3333333333, .7698003590, -.5443310539);
\draw (-.3333333333, .7698003590, -.5443310539)--(.4444444444, 1.283000598, -.9072184234);
\tikzstyle{every node}=[draw=black,fill=yellow!40!white,thick,
  shape=circle,minimum size=1em,inner sep=0.4];	
\node at (0., 0., 0.) {1};
\node at (1., 0., 0.) {2};
\node at (.5000000000, -.6735753142, -.5443310539) {3};
\node at (.5000000000, .2886751347, .8164965809) {4};
\node at (.5000000000, .2886751347, -.8164965809) {5};
\node at (.5000000000, .8660254040, 0.) {6};
\node at (1.333333333, .7698003590, -.5443310539) {7};
\node at (-.3888888889, -.1603750748, -.9072184234) {8};
\node at (-.3888888889, .8018753741, .4536092117) {9};
\node at (-.3333333333, .7698003590, -.5443310539) {\footnotesize {10}};
\node at (.4444444444, 1.283000598, -.9072184234) {\footnotesize {11}};
\end{tikzpicture}
\begin{tikzpicture}
\tikzset{
    table/.style={
        matrix of nodes,
        row sep=-\pgflinewidth,
        column sep=-\pgflinewidth,
        nodes={
            rectangle,
            draw=black,
            align=center
        },
        minimum height=2.5em,
        text depth=1ex,
        text height=2ex,
        nodes in empty cells,
        every even row/.style={
            nodes={fill=gray!20}
        },
        column 1/.style={
            nodes={text width=2em,font=\bfseries}
        },
        row 1/.style={
            nodes={
                fill=black,
                text=white,
                font=\bfseries
            }
        }
    }
}

\matrix (first) [table,text width=7em]
{
& X & Y & Z\\
1 & $0$ & $0$ & $0$ \\
2 & $1$ & $0$ & $0$ \\
3 & $1/2$ & $-7\sqrt{3}/18$ & $-2\sqrt{6}/9$ \\
4 & $1/2$ & $\sqrt{3}/6$ & $\sqrt{6}/3$ \\
5 & $1/2$ & $\sqrt{3}/6$ & $-\sqrt{6}/3$ \\
6 & $1/2$ & $\sqrt{3}/2$ & $0$ \\
7 & $4/3$ & $4\sqrt{3}/9$ & $-2\sqrt{6}/9$ \\
8 & $-7/18$ & $-5\sqrt{3}/54$ & $-10\sqrt{6}/27$ \\
9 & $-7/18$ & $25\sqrt{3}/54$ & $5\sqrt{6}/27$ \\
10 & $-1/3$ & $4\sqrt{3}/9$ & $-2\sqrt{6}/9$ \\
11 & $4/9$ & $20\sqrt{3}/27$ &$-10\sqrt{6}/27$ \\
};
\end{tikzpicture}

\begin{figure}[!htb]
\centering
\caption{Goldner-Harary Graph}
\label{}
\end{figure}

This was probably the easiest graph of all since it had some many regular tetrahedra there was basically no degree of freedom.
We suspect the vertices of this graph represent the vertices of some polyhedron; unfortunately, we could not get a hold
of Goldner and Harary's article since it appeared in some obscure Malaysian journal.
\newpage
\hspace{-2cm}
\tdplotsetmaincoords{58}{159}
\begin{tikzpicture}[tdplot_main_coords,
            cube/.style={very thick,black},
			grid/.style={gray,very thin,opacity=0.05},
			axis/.style={->,blue,very thick},
            line width=1.0pt,scale=3.8]
\foreach \x in {-0.8,-0.4,...,1.2}
\foreach \y in {-0.8,-0.4,...,1.2}
\foreach \z in {-0.8,-0.4,...,1.2}
	{
			\draw[grid] (\x,-0.8,-0.8) -- (\x,1.2,-0.8);
			\draw[grid] (-0.8,\y,-0.8) -- (1.2,\y,-0.8);
            \draw[grid] (-0.8,-0.8,\z) -- (1.2,-0.8,\z);
            \draw[grid] (-0.8,-0.8,\z) -- (-0.8,1.2,\z);
           \draw[grid] (\x,-0.8,-0.8) -- (\x,-0.8,1.2);
            \draw[grid] (-0.8,\y,-0.8) -- (-0.8,\y,1.2);
		}
\draw[axis] (0,0,0) -- (1,0,0) node[anchor=north]{$x$};
	\draw[axis] (0,0,0) -- (0,1.5,0) node[anchor=west]{$y$};
	\draw[axis] (0,0,0) -- (0,0,1.3) node[anchor=east]{$z$};
    \draw (0,0,0)--(-.3600000000,.4800000000,.8000000000);
    \draw (0,0,0)--(.3600000000, .4800000000, -.8000000000);
    \draw (0,0,0)--(.3600000000, -.4800000000, -.8000000000);
    \draw (-.3600000000,.4800000000,.8000000000)--(0., .9600000000, 0.);
    \draw (-.3600000000,.4800000000,.8000000000)--(-.7200000000, 0., 0);
    \draw (0., .9600000000, 0.)--(.3600000000, .4800000000, -.8000000000);
    \draw (0., .9600000000, 0.)--(-.2800000000, 0., 0.);
    \draw (0., .9600000000, 0.)--(.2800000000, 0., 0);
    \draw (.3600000000, .4800000000, -.8000000000)--(.7200000000, 0., 0.);
    \draw (.7200000000, 0., 0.)--(-.2800000000, 0., 0.);
    \draw (.7200000000, 0., 0.)--(.3600000000, -.4800000000, -.8000000000);
    \draw (-.2800000000, 0., 0.)--(0., -.9600000000, 0.);
    \draw (.3600000000, -.4800000000, -.8000000000)--(0., -.9600000000, 0.);
    \draw (0., -.9600000000, 0.)--(-.3600000000, -.4800000000, .8000000000);
    \draw (0., -.9600000000, 0.)--(.2800000000, 0., 0);
    \draw (0,0,0)--(-.3600000000, -.4800000000, .8000000000);
    \draw (-.3600000000, -.4800000000, .8000000000)--(-.7200000000, 0., 0);
    \draw (-.7200000000, 0., 0)--(.2800000000, 0., 0);
\tikzstyle{every node}=[draw=black,fill=yellow!40!white,thick,
  shape=circle,minimum size=1.2em,inner sep=0.4];	
    \node at (0,0,0) {1};
    \node at (-.3600000000,.4800000000,.8000000000) {2};
    \node at (0., .9600000000, 0.) {3};
    \node at (.3600000000, .4800000000, -.8000000000) {4};
    \node at (.7200000000, 0., 0.) {5};
    \node at (-.2800000000, 0., 0.) {6};
    \node at (.3600000000, -.4800000000, -.8000000000) {7} ;
    \node at (0., -.9600000000, 0.) {8};
    \node at (-.3600000000, -.4800000000, .8000000000) {9};
    \node at (-.7200000000, 0., 0) {10};
    \node at (.2800000000, 0., 0) {11};
\end{tikzpicture}
\begin{tikzpicture}
\tikzset{
    table/.style={
        matrix of nodes,
        row sep=-\pgflinewidth,
        column sep=-\pgflinewidth,
        nodes={
            rectangle,
            draw=black,
            align=center
        },
        minimum height=2.0em,
        text depth=1ex,
        text height=2ex,
        nodes in empty cells,
        every even row/.style={
            nodes={fill=gray!20}
        },
        column 1/.style={
            nodes={text width=3em,font=\bfseries}
        },
        row 1/.style={
            nodes={
                fill=black,
                text=white,
                font=\bfseries
            }
        }
    }
}

\matrix (first) [table,text width=6em]
{
& X & Y & Z \\
1 & $0$ & $0$ & $0$ \\
2 & $-9/25$ & $12/25$ & $4/5$ \\
3 & $0$ & $24/25$ & $0$ \\
4 & $9/25$ & $12/25$ & $-4/5$ \\
5 & $18/25$ & $0$ & $0$ \\
6 & $-7/25$ & $0$ & $0$ \\
7 & $9/25$ & $-12/25$ & $-4/5$ \\
8 & $0$ & $-24/25$ & $0$ \\
9 & $-9/25$ & $-12/25$ & $4/5$ \\
10 & $-18/25$ & $0$ & $0$ \\
11 & $7/25$ & $0$ & $0$ \\
};
\end{tikzpicture}

\begin{figure}[!htb]
\centering
\caption{Herschel Graph}
\label{}
\end{figure}
This is probably one of the prettiest embeddings if one judges by the coordinates.
Place $A_1$ at the origin and then use the central symmetry of the graph to define $A_{5+k}=-A_k$ for all $2\le k\le 6$.
Further, assume that the points $A_3$, $A_5$, and $A_6$ (and therefore $A_8$, $A_{10}$, and $A_{11}$) lie in the plane $z=0$.
Finally, add the restriction that the points $A_1,\,A_2,\,A_3,\,A_4$ are actually coplanar, and therefore form a rhombus.
This leads to a system of $8$ equations and $8$ unknowns whose solution leads to the embedding presented above.
\newpage
\hspace{-2cm}
\tdplotsetmaincoords{58}{157}
\begin{tikzpicture}[tdplot_main_coords,
            cube/.style={very thick,black},
			grid/.style={gray,very thin,opacity=0.05},
			axis/.style={->,blue,very thick},
    line width=1.0pt,scale=3.5]
 \draw[axis] (0,0,0) -- (1.5,0,0) node[anchor=west]{$x$};
	\draw[axis] (0,0,0) -- (0,1.8,0) node[anchor=west]{$y$};
	\draw[axis] (0,0,0) -- (0,0,1.5) node[anchor=west]{$z$};
\foreach \x in {-0.8,-0.4,...,1.2}
\foreach \y in {-0.8,-0.4,...,1.2}
\foreach \z in {-0.8,-0.4,...,1.2}
	{
			\draw[grid] (\x,-0.8,-0.8) -- (\x,1.2,-0.8);
			\draw[grid] (-0.8,\y,-0.8) -- (1.2,\y,-0.8);
            \draw[grid] (-0.8,-0.8,\z) -- (1.2,-0.8,\z);
            \draw[grid] (-0.8,-0.8,\z) -- (-0.8,1.2,\z);
           \draw[grid] (\x,-0.8,-0.8) -- (\x,-0.8,1.2);
            \draw[grid] (-0.8,\y,-0.8) -- (-0.8,\y,1.2);
		}
\draw (.5773502693, 0., .5000000000)--(-.2886751347, .5000000000, .5000000000);
\draw (.5773502693, 0., .5000000000)--(-.2886751347, -.5000000000, .5000000000);
\draw (.5773502693, 0., .5000000000)--(.5773502693, 0., -.5000000000);
\draw (.5773502693, 0., .5000000000)--(.4978909578, .8623724358, 0.);
\draw (.5773502693, 0., .5000000000)--(.4978909578, -.8623724358, 0.);
\draw (-.2886751347, .5000000000, .5000000000)--(-.2886751347, -.5000000000, .5000000000);
\draw (-.2886751347, .5000000000, .5000000000)--(-.9957819157, 0., 0.);
\draw (-.2886751347, .5000000000, .5000000000)--(-.2886751347, .5000000000, -.5000000000);
\draw (-.2886751347, .5000000000, .5000000000)--(.4978909578, .8623724358, 0.);
\draw (-.2886751347, -.5000000000, .5000000000)--(-.9957819157, 0., 0.);
\draw (-.2886751347, -.5000000000, .5000000000)--(-.2886751347, -.5000000000, -.5000000000);
\draw (-.2886751347, -.5000000000, .5000000000)--(.4978909578, -.8623724358, 0.);
\draw (-.9957819157, 0., 0.)--(-.2886751347, -.5000000000, -.5000000000);
\draw (-.9957819157, 0., 0.)--(-.2886751347, .5000000000, -.5000000000);
\draw (-.2886751347, -.5000000000, -.5000000000)--(-.2886751347, .5000000000, -.5000000000);
\draw (-.2886751347, -.5000000000, -.5000000000)--(.5773502693, 0., -.5000000000);
\draw (-.2886751347, -.5000000000, -.5000000000)--(.4978909578, -.8623724358, 0.);
\draw (-.2886751347, .5000000000, -.5000000000)--(.5773502693, 0., -.5000000000);
\draw (-.2886751347, .5000000000, -.5000000000)--(.4978909578, .8623724358, 0.);
\draw (.5773502693, 0., -.5000000000)--(.4978909578, .8623724358, 0.);
\draw (.5773502693, 0., -.5000000000)--(.4978909578, -.8623724358, 0.);
\tikzstyle{every node}=[draw=black,fill=yellow!40!white,thick,
  shape=circle,minimum size=1em,inner sep=0.4];	
\node at (.5773502693, 0., .5000000000) {1};
\node at (-.2886751347, .5000000000, .5000000000) {2};
\node at (-.2886751347, -.5000000000, .5000000000) {3};
\node at (-.9957819157, 0., 0.) {4};
\node at (-.2886751347, -.5000000000, -.5000000000) {5};
\node at (-.2886751347, .5000000000, -.5000000000) {6};
\node at (.5773502693, 0., -.5000000000) {7};
\node at (.4978909578, .8623724358, 0.) {8};
\node at (.4978909578, -.8623724358, 0.) {9};
\end{tikzpicture}
\begin{tikzpicture}
\tikzset{
    table/.style={
        matrix of nodes,
        row sep=-\pgflinewidth,
        column sep=-\pgflinewidth,
        nodes={
            rectangle,
            draw=black,
            align=center
        },
        minimum height=2.0em,
        text depth=1ex,
        text height=2ex,
        nodes in empty cells,
        every even row/.style={
            nodes={fill=gray!20}
        },
        column 1/.style={
            nodes={text width=3em,font=\bfseries}
        },
        row 1/.style={
            nodes={
                fill=black,
                text=white,
                font=\bfseries
            }
        }
    }
}

\matrix (first) [table,text width=6em]
{
& X & Y & Z \\
1 & $\sqrt{3}/3$ & $0$ & $1/2$ \\
2 & $-\sqrt{3}/6$ & $1/2$ & $1/2$ \\
3 & $-\sqrt{3}/6$ & $-1/2$ & $1/2$ \\
4 & {\small $-\sqrt{3}/6-\sqrt{2}/2$} & $0$ & $0$ \\
5 & $-\sqrt{3}/6$ & $-1/2$ & $-1/2$ \\
6 & $-\sqrt{3}/6$ & $1/2$ & $-1/2$ \\
7 & $\sqrt{3}/3$ & $0$ & $-1/2$ \\
8 & $\sqrt{3}/12+\sqrt{2}/4$ & $1/4+\sqrt{6}/4$ & $0$ \\
9 & $\sqrt{3}/12+\sqrt{2}/4$ & $-1/4-\sqrt{6}/6$ & $0$ \\
};
\end{tikzpicture}

\begin{figure}[!htb]
\centering
\caption{Fritsch Graph}
\label{}
\end{figure}

For this graph we employed the rigidity of the graph and a clever observation to obtain the embedding.
We postulated that because of the rigidity of the graph, the Euclidean embedding would result in a
triangular prism attached with $3$ regular pyramids.
We took $A_5$, $A_6$, $A_7$ as the base of our prism. Since this would form a equilateral triangle, we locate
$A_5=[-\sqrt{3}/6, -1/2, -1/2], A_6=[-\sqrt{3}/6, 1/2, -1/2], A_7=[\sqrt{3}/3, 0, -1/2]$.
We then selected $A_1, A_2, A_3$ as the opposite vertices of the prism. Since $|A_1-A_7|=1, |A_2-A_6|=1, |A_3-A_5|=1$, we selected
$A_1=[-\sqrt{3}/6, -1/2, 1/2], A_2=[-\sqrt{3}/6, 1/2, 1/2], A_3=[\sqrt{3}/3, 0, 1/2]$.
Since $A_4$ is part of the regular pyramid consisting of $A_2, A_3, A_5, A_6$, we are able to easily compute as $A_4=[-\sqrt{3}/6-\sqrt{2}/2,0]$. Vertices $A_8$ and $A_9$ are $60^{\circ}$ rotations of $A_4$.
\newpage

\hspace{-2.3cm}
\tdplotsetmaincoords{60}{125}
\begin{tikzpicture}[tdplot_main_coords,
            cube/.style={very thick,black},
			grid/.style={gray,very thin,opacity=0.05},
			axis/.style={->,blue,very thick},
    line width=1.0pt,scale=3.5]
    \draw[axis] (0,0,0) -- (2,0,0) node[anchor=west]{$x$};
	\draw[axis] (0,0,0) -- (0,1.2,0) node[anchor=west]{$y$};
	\draw[axis] (0,0,0) -- (0,0,2) node[anchor=west]{$z$};
\foreach \x in {-0.8,-0.4,...,1.2}
\foreach \y in {-0.8,-0.4,...,1.2}
\foreach \z in {-0.8,-0.4,...,1.2}
	{
			\draw[grid] (\x,-0.8,-0.8) -- (\x,1.2,-0.8);
			\draw[grid] (-0.8,\y,-0.8) -- (1.2,\y,-0.8);
            \draw[grid] (-0.8,-0.8,\z) -- (1.2,-0.8,\z);
            \draw[grid] (-0.8,-0.8,\z) -- (-0.8,1.2,\z);
           \draw[grid] (\x,-0.8,-0.8) -- (\x,-0.8,1.2);
            \draw[grid] (-0.8,\y,-0.8) -- (-0.8,\y,1.2);
		}
\draw (0., 0., 0.)--(-.1004057078, .3090169946, .9457416090);
\draw (0., 0., 0.)--(-.1004057078, -.3090169946, .9457416090);
\draw (0., 0., 0.)--(.2628655562, .1909830058, .9457416090);
\draw (0., 0., 0.)--(-.3249196964, 0., .9457416090);
\draw (0., 0., 0.)--(.2628655562, -.1909830058, .9457416090);
\draw (-.1004057078, .3090169946, .9457416090)--(.8506508080, 0., .9457416082);
\draw (-.1004057078, .3090169946, .9457416090)--(-.6881909600, -.5000000000, .9457416082);
\draw (-.1004057078, -.3090169946, .9457416090)--(.8506508080, 0., .9457416082);
\draw (-.1004057078, -.3090169946, .9457416090)--(-.6881909600, .5000000000, .9457416082);
\draw (.2628655562, .1909830058, .9457416090)--(-.6881909600, .5000000000, .9457416082);
\draw (.2628655562, .1909830058, .9457416090)--(.2628655554, -.8090169940, .9457416082);
\draw (-.3249196964, 0., .9457416090)--(.2628655554, -.8090169940, .9457416082);
\draw (-.3249196964, 0., .9457416090)--(-.6881909600, -.5000000000, .9457416082);
\draw (.2628655562, -.1909830058, .9457416090)--(-.6881909600, -.5000000000, .9457416082);
\draw (.2628655562, -.1909830058, .9457416090)--(.2628655554, .8090169940, .9457416082);
\draw (.8506508080, 0., .9457416082)--(.2628655554, -.8090169940, .9457416082);
\draw (.8506508080, 0., .9457416082)--(.2628655554, .8090169940, .9457416082);
\draw (-.6881909600, .5000000000, .9457416082)--(.2628655554, .8090169940, .9457416082);
\draw (-.6881909600, .5000000000, .9457416082)--(-.6881909600, -.5000000000, .9457416082);
\draw (.2628655554, -.8090169940, .9457416082)--(-.6881909600, -.5000000000, .9457416082);
\draw (-.3249196964, 0., .9457416090)--(.2628655554, .8090169940, .9457416082);
\tikzstyle{every node}=[draw=black,fill=yellow!50!white,thick,
  shape=circle,minimum size=1.3em,inner sep=0.4];	
\node at (0., 0., 0.) {1};
\node at (-.1004057078, .3090169946, .9457416090) {2};
\node at (-.1004057078, -.3090169946, .9457416090) {3};
\node at (.2628655562, .1909830058, .9457416090) {4};
\node at (-.3249196964, 0., .9457416090) {5};
\node at (.2628655562, -.1909830058, .9457416090) {6};
\node at (.8506508080, 0., .9457416082) {7};
\node at (-.6881909600, .5000000000, .9457416082) {8};
\node at (.2628655554, -.8090169940, .9457416082) {9};
\node at (.2628655554, .8090169940, .9457416082) {10};
\node at (-.6881909600, -.5000000000, .9457416082) {11};
\end{tikzpicture}
\begin{tikzpicture}
\tikzset{
    table/.style={
        matrix of nodes,
        row sep=-\pgflinewidth,
        column sep=-\pgflinewidth,
        nodes={
            rectangle,
            draw=black,
            align=center
        },
        minimum height=2.0em,
        text depth=1.5ex,
        text height=2.7ex,
        nodes in empty cells,
        every even row/.style={
            nodes={fill=gray!20}
        },
        column 1/.style={
            nodes={text width=1em,font=\bfseries}
        },
        row 1/.style={
            nodes={
                fill=black,
                text=white,
                font=\bfseries
            }
        }
    }
}

\matrix (first) [table,text width=6em]
{
& X & Y & Z \\
1 & $0$ & $0$ & $0$ \\
2 & {\small $-2b\cos(2\pi/5)$} & {\small $2b\sin(2\pi/5)$} & $c$ \\
3 & {\small $-2b\cos(2\pi/5)$} & {\small $-2b\sin(2\pi/5)$} & $c$ \\
4 & $2b\cos(2\pi/5)$ & $2b\sin(2\pi/5)$ & $c$ \\
5 & $-2b$ & $0$ & $c$ \\
6 & $2b\cos(\pi/5)$ & $2b\sin(\pi/5)$ & $c$ \\
7 & $1/2\sin(\pi/5)$ & $0$ & $c$ \\
8 & $-\cot(\pi/5)/2$ & $1/2$ & $c$ \\
9 & {\tiny $\cot(2\pi/5)\cos(\pi/5)$} & $-\cos(\pi/5)$ & $c$ \\
10 & {\tiny $\cot(2\pi/5)\cos(\pi/5)$} & $\cos(\pi/5)$ & $c$ \\
11 & $-\cot(\pi/5)/2$ & $-1/2$ & $c$ \\
};

\end{tikzpicture}

\begin{figure}[!htb]
\centering
\caption{Gr\"{o}tzsch Graph: $b=1/2\sqrt{5+2\sqrt{5}}$, $c=\sqrt{5+3\sqrt{5}}(5-\sqrt{5})/10$}
\label{}
\end{figure}

Very easy embedding using the $5$-fold symmetry. Set vertex $1$ at the origin; the remaining vertices
form $2$ coplanar concentric regular pentagons.
Let $A_5$ be in the $xz$ plane.  The coordinates of the points $A_3$, $A_6$, $A_4$, and $A_2$, are then obtained by $72^{\circ}$ rotations. Let $A_7$ have the same $z$ coordinate as $A_5$. So far we have four variables, two from $A_5$ and two from $A_7$.
The points $A_{10}$, $A_8$, $A_{11}$, and $A_9$, are then obtained by rotating $A_7$ in $72^{\circ}$ increments.
The number of equations is also four as we have four types of edges: $\|A_1-A_5\|=1$, $\|A_5-A_9\|=1$, $\|A_5-A_{10}\|=1$ and $\|A_7-A_{10}\|=1$. Solving this system produces the embedding shown above.
\newpage
\hspace{-2.5cm}
\tdplotsetmaincoords{40}{40}
\begin{tikzpicture}[tdplot_main_coords,
            cube/.style={very thick,black},
			grid/.style={gray,very thin,opacity=0.05},
			axis/.style={->,blue,very thick},
    line width=1.0pt,scale=3.7]
 \draw[axis] (0,0,0) -- (1.2,0,0) node[anchor=west]{$x$};
	\draw[axis] (0,0,0) -- (0,3.0,0) node[anchor=west]{$y$};
	\draw[axis] (0,0,0) -- (0,0,1.5) node[anchor=west]{$z$};
\foreach \x in {-0.2,0,...,1.2}
\foreach \y in {-0.2,0,...,1.2}
\foreach \z in {-0.2,0,...,1.2}
	{
			\draw[grid] (\x,-0.2,-0.2) -- (\x,1.2,-0.2);
			\draw[grid] (-0.2,\y,-0.2) -- (1.2,\y,-0.2);
            \draw[grid] (-0.2,-0.2,\z) -- (1.2,-0.2,\z);
            \draw[grid] (-0.2,-0.2,\z) -- (-0.2,1.2,\z);
           \draw[grid] (\x,-0.2,-0.2) -- (\x,-0.2,1.2);
            \draw[grid] (-0.2,\y,-0.2) -- (-0.2,\y,1.2);
		}
\draw (0., 2.800000000, 0.)--(-.7783817294, 2.200000000, .1847210963);
\draw (0., 2.800000000, 0.)--(.7783817294, 2.200000000, .1847210963);
\draw (0., 2.800000000, 0.)--(.6000000000, 2., 0.);
\draw (0., 2.800000000, 0.)--(-.6000000000, 2., 0.);
\draw (-.7783817294, 2.200000000, .1847210963)--(-.1783817294, 1.400000000, .1847210963);
\draw (-.7783817294, 2.200000000, .1847210963)--(0., 1.600000000, 0.);
\draw (-.7783817294, 2.200000000, .1847210963)--(-.1186663608, 1.518666361, .5018332290);
\draw (-.1783817294, 1.400000000, .1847210963)--(-.1574374633, .6000000000, .7843554329);
\draw (-.1783817294, 1.400000000, .1847210963)--(.6000000000, 2., 0.);
\draw (-.1783817294, 1.400000000, .1847210963)--(.6000000000, .8000000000, 0.);
\draw (-.1574374633, .6000000000, .7843554329)--(0., 0., 0.);
\draw (-.1574374633, .6000000000, .7843554329)--(0., 1.200000000, 0.);
\draw (-.1574374633, .6000000000, .7843554329)--(.1186663608, 1.518666361, .50183322);
\draw (0., 0., 0.)--(.1574374633, .6000000000, .7843554329);
\draw (0., 0., 0.)--(.6000000000, .8000000000, 0.);
\draw (0., 0., 0.)--(-.6000000000, .8000000000, 0.);
\draw (.1574374633, .6000000000, .7843554329)--(.1783817294, 1.400000000, .1847210963);
\draw (.1574374633, .6000000000, .7843554329)--(0., 1.200000000, 0.);
\draw (.1574374633, .6000000000, .7843554329)--(-.1186663608, 1.518666361, .5018332290);
\draw (.1783817294, 1.400000000, .1847210963)--(.7783817294, 2.200000000, .1847210963);
\draw (.1783817294, 1.400000000, .1847210963)--(-.6000000000, .8000000000, 0.);
\draw (.1783817294, 1.400000000, .1847210963)--(-.6000000000, 2., 0.);
\draw (.7783817294, 2.200000000, .1847210963)--(0., 1.600000000, 0.);
\draw (.7783817294, 2.200000000, .1847210963)--(.1186663608, 1.518666361, .50183322);
\draw (.6000000000, 2., 0.)--(0., 1.200000000, 0.);
\draw (.6000000000, 2., 0.)--(-.1186663608, 1.518666361, .5018332290);
\draw (.6000000000, .8000000000, 0.)--(0., 1.600000000, 0.);
\draw (.6000000000, .8000000000, 0.)--(.1186663608, 1.518666361, .5018332290);
\draw (-.6000000000, .8000000000, 0.)--(0., 1.600000000, 0.);
\draw (-.6000000000, .8000000000, 0.)--(-.1186663608, 1.518666361, .5018332290);
\draw (-.6000000000, 2., 0.)--(0., 1.200000000, 0.);
\draw (-.6000000000, 2., 0.)--(.1186663608, 1.518666361, .5018332290);
\node at (0.1, 2.800000000, 0.) {1};
\node at (-.8783817294, 2.200000000, .1847210963) {2};
\node at (-.1783817294, 1.340000000, .1047210963) {3};
\node at (-.1574374633, .6000000000, .8843554329) {4};
\node at (0.05, -0.05, -0.05) {5};
\node at (.1574374633, .6000000000, .6543554329) {6};
\node at (.2783817294, 1.400000000, .1847210963) {7};
\node at (.8583817294, 2.150000000, .1847210963) {8};
\node at (.6500000000, 1.95, -0.05) {9};
\node at (.6500000000, .7500000000, -0.05) {10};
\node at (-.6000000000, 0.9000000000, -0.17) {11};
\node at (-.7500000000, 2., 0.) {12};
\node at (0.0, 1.6500000000, 0.1) {13};
\node at (0.10, 1.120000000, 0.) {14};
\node at (-.1186663608, 1.478666361, .4118332290) {15};
\node at (.1186663608, 1.618666361, .50183322) {16};
\tikzstyle{every node}=[draw=red,fill=red!20!,shape=circle,scale=0.2]
\node at (0., 2.800000000, 0.) {};
\node at (-.7783817294, 2.200000000, .1847210963) {};
\node at (-.1783817294, 1.400000000, .1847210963) {};
\node at (-.1574374633, .6000000000, .7843554329) {};
\node at (0., 0., 0.) {};
\node at (.1574374633, .6000000000, .7843554329) {};
\node at (.1783817294, 1.400000000, .1847210963) {};
\node at (.7783817294, 2.200000000, .1847210963) {};
\node at (.6000000000, 2., 0.) {};
\node at (.6000000000, .8000000000, 0.) {};
\node at (-.6000000000, .8000000000, 0.) {};
\node at (-.6000000000, 2., 0.) {};
\node at (0., 1.600000000, 0.) {};
\node at (0., 1.200000000, 0.) {};
\node at (-.1186663608, 1.518666361, .5018332290) {};
\node at (.1186663608, 1.518666361, .50183322) {};
\end{tikzpicture}
\begin{tikzpicture}
\tikzset{
    table/.style={
        matrix of nodes,
        row sep=-\pgflinewidth,
        column sep=-\pgflinewidth,
        nodes={
            rectangle,
            draw=black,
            align=center
        },
        minimum height=2.0em,
        text depth=1ex,
        text height=1ex,
        nodes in empty cells,
        every even row/.style={
            nodes={fill=gray!20}
        },
        column 1/.style={
            nodes={text width=1em,font=\bfseries}
        },
        row 1/.style={
            nodes={
                fill=black,
                text=white,
                font=\bfseries
            }
        }
    }
}

\matrix (first) [table,text width=8.0em]
{
& X & Y & Z \\
1 & $0$ & $14/5$ & $0$ \\
2 & {\scriptsize $-4(s^{2}-1)/5(s^{2}+1)$} & $11/5$ & $8s/5(s^{2}+1)$ \\
3 & {\scriptsize $-(s^{2}-7)/(s^{2}+1)$} & $7/5$ & $8s/5(s^{2}+1)$ \\
4 &{\scriptsize $-4(s^{2}-49)/5(s^{2}+49)$}& $3/5$ & {\footnotesize $56s/5(s^{2}+49)$} \\
5 & $0$ & $0$ & $0$ \\
6 & {\scriptsize $4(s^{2}-49)/5(s^{2}+49)$} & $3/5$ & {\footnotesize $56s/5(s^{2}+49)$} \\
7 & {\scriptsize $(s^{2}-7)/5(s^2+1)$} & $7/5$ & $8s/5(s^{2}+1)$ \\
8 & {\scriptsize $4(s^{2}+1)/5(s^{2}+1)$} & $11/5$ & $8s/5(s^{2}+1)$ \\
9 & $3/5$ & $2$ & $0$ \\
10 & $3/5$ & $4/5$ & $0$ \\
11 & $-3/5$ & $4/5$ & $0$ \\
12 & $-3/5$ & $2$ & $0$ \\
13 & $0$ & $8/5$ & $0$ \\
14 & $0$ & $6/5$ & $0$ \\
15 & {\scriptsize $-21(s^{2}-7)/160s^{2}$} & {\tiny $49(5s^{2}-3)/160s^2$} & {\scriptsize $7(s^{2}+25)/160s$} \\
16 & {\scriptsize $21(s^{2}-7)/160s^{2}$} & {\tiny $49(5s^{2}-3)/160s^{2}$} & {\scriptsize $7(s^{2}+25)/160s$} \\
};
\end{tikzpicture}

\begin{figure}[!htb]
\centering
\caption{Hoffman Graph: $s$=$\sqrt{6986+14\sqrt{273697}}/14$}
\label{}
\end{figure}

Use the symmetry of the graph.
Start by selecting vertices $1$, $5$, $13$, and $14$ on the $y$-axis and the remaining points symmetric in pairs
with respect to the $yz$-plane: $3\leftrightarrow 7$, $4\leftrightarrow 6$, $2\leftrightarrow 8$, $9\leftrightarrow 12$, $10\leftrightarrow 11$, $15\leftrightarrow 16$.
\newpage
\centering
\tdplotsetmaincoords{58}{130}
\begin{tikzpicture}[tdplot_main_coords,
            cube/.style={very thick,black},
			grid/.style={gray,very thin,opacity=0.05},
			axis/.style={->,blue,very thick},
    line width=1.0pt,scale=3.3]
 \draw[axis] (0,0,0) -- (1.5,0,0) node[anchor=west]{$x$};
	\draw[axis] (0,0,0) -- (0,1.8,0) node[anchor=west]{$y$};
	\draw[axis] (0,0,0) -- (0,0,0.5) node[anchor=west]{$z$};
\foreach \x in {-0.8,-0.4,...,1.2}
\foreach \y in {-0.8,-0.4,...,1.2}
\foreach \z in {-0.8,-0.4,...,0.4}
	{
			\draw[grid] (\x,-0.8,-0.8) -- (\x,1.2,-0.8);
			\draw[grid] (-0.8,\y,-0.8) -- (1.2,\y,-0.8);
            \draw[grid] (-0.8,-0.8,\z) -- (1.2,-0.8,\z);
            \draw[grid] (-0.8,-0.8,\z) -- (-0.8,1.2,\z);
           \draw[grid] (\x,-0.8,-0.8) -- (\x,-0.8,0.4);
            \draw[grid] (-0.8,\y,-0.8) -- (-0.8,\y,0.4);
		}
\draw (0., 0., 0.)--(1., 0., 0.);
\draw (0., 0., 0.)--(.50000000000000000000000000000000000000000000000000, .86602540378443864676372317075293618347140262690520, 0.);
\draw (0., 0., 0.)--(.50000000000000000000000000000000000000000000000000, .70996725453148610937235272105839245675313139363253, -.49592993203982356170374361442954598041731399972938);
\draw (0., 0., 0.)--(-.46376374908530441738621882110796960007191521872411, .84510439456441632930590039236874174086997878496293, -.26593560747341475551243840504831172694944221365704);
\draw (1., 0., 0.)--(.50000000000000000000000000000000000000000000000000, .86602540378443864676372317075293618347140262690520, 0.);
\draw (1., 0., 0.)--(.50000000000000000000000000000000000000000000000000, .70996725453148610937235272105839245675313139363253, -.49592993203982356170374361442954598041731399972938);
\draw (1., 0., 0.)--(.71093242412093931363351251934884378021978648378953, .41045702644185750569725934296678974357362854919682, -.86484967827935958644932861307814066351067871619060);
\draw (1., 0., 0.)--(0.88984975091071642355704808526684294906800983992502e-1, -.19218828437142689074155991923309588879620927545075, -.36484967827935958644932861307814066351067871619058);
\draw (.50000000000000000000000000000000000000000000000000, .86602540378443864676372317075293618347140262690520, 0.)--(-.46376374908530441738621882110796960007191521872411, .84510439456441632930590039236874174086997878496293, -.26593560747341475551243840504831172694944221365704);
\draw (.50000000000000000000000000000000000000000000000000, .86602540378443864676372317075293618347140262690520, 0.)--(.71093242412093931363351251934884378021978648378953, .41045702644185750569725934296678974357362854919682, -.86484967827935958644932861307814066351067871619060);
\draw (.50000000000000000000000000000000000000000000000000, .86602540378443864676372317075293618347140262690520, 0.)--(.21093242412093931363351251934884378021978648378953, 1.2764824302262961524609825137197259270450311761020, -.86484967827935958644932861307814066351067871619060);
\draw (.50000000000000000000000000000000000000000000000000, .70996725453148610937235272105839245675313139363253, -.49592993203982356170374361442954598041731399972938)--(-.46376374908530441738621882110796960007191521872411, .84510439456441632930590039236874174086997878496293, -.26593560747341475551243840504831172694944221365704);
\draw (.50000000000000000000000000000000000000000000000000, .70996725453148610937235272105839245675313139363253, -.49592993203982356170374361442954598041731399972938)--(0.88984975091071642355704808526684294906800983992502e-1, -.19218828437142689074155991923309588879620927545075, -.36484967827935958644932861307814066351067871619058);
\draw (.50000000000000000000000000000000000000000000000000, .70996725453148610937235272105839245675313139363253, -.49592993203982356170374361442954598041731399972938)--(-.26811812545734779130689058944601519996404239063795, .42255219728220816465295019618437087043498939248147, -1.0681074658615250936386579889886538538269692085900);
\draw (-.46376374908530441738621882110796960007191521872411, .84510439456441632930590039236874174086997878496293, -.26593560747341475551243840504831172694944221365704)--(.21093242412093931363351251934884378021978648378953, 1.2764824302262961524609825137197259270450311761020, -.86484967827935958644932861307814066351067871619060);
\draw (.71093242412093931363351251934884378021978648378953, .41045702644185750569725934296678974357362854919682, -.86484967827935958644932861307814066351067871619060)--(0.88984975091071642355704808526684294906800983992502e-1, -.19218828437142689074155991923309588879620927545075, -.36484967827935958644932861307814066351067871619058);
\draw (.71093242412093931363351251934884378021978648378953, .41045702644185750569725934296678974357362854919682, -.86484967827935958644932861307814066351067871619060)--(.21093242412093931363351251934884378021978648378953, 1.2764824302262961524609825137197259270450311761020, -.86484967827935958644932861307814066351067871619060);
\draw (.71093242412093931363351251934884378021978648378953, .41045702644185750569725934296678974357362854919682, -.86484967827935958644932861307814066351067871619060)--(-.26811812545734779130689058944601519996404239063795, .42255219728220816465295019618437087043498939248147, -1.0681074658615250936386579889886538538269692085900);
\draw (0.88984975091071642355704808526684294906800983992502e-1, -.19218828437142689074155991923309588879620927545075, -.36484967827935958644932861307814066351067871619058)--(-.26811812545734779130689058944601519996404239063795, .42255219728220816465295019618437087043498939248147, -1.0681074658615250936386579889886538538269692085900);
\draw (.21093242412093931363351251934884378021978648378953, 1.2764824302262961524609825137197259270450311761020, -.86484967827935958644932861307814066351067871619060)--(-.26811812545734779130689058944601519996404239063795, .42255219728220816465295019618437087043498939248147, -1.0681074658615250936386579889886538538269692085900);
\tikzstyle{every node}=[draw=black,fill=yellow!50!white,thick,
  shape=circle,minimum size=1em,inner sep=0.2];	
\node at (0., 0., 0.) {1};
\node at (1., 0., 0.) {2};
\node at (.50000000000000000000000000000000000000000000000000, .86602540378443864676372317075293618347140262690520, 0.) {3};
\node at (.50000000000000000000000000000000000000000000000000, .70996725453148610937235272105839245675313139363253, -.49592993203982356170374361442954598041731399972938) {4};
\node at (-.46376374908530441738621882110796960007191521872411, .84510439456441632930590039236874174086997878496293, -.26593560747341475551243840504831172694944221365704) {5};
\node at (.71093242412093931363351251934884378021978648378953, .41045702644185750569725934296678974357362854919682, -.86484967827935958644932861307814066351067871619060) {6};
\node at (0.88984975091071642355704808526684294906800983992502e-1, -.19218828437142689074155991923309588879620927545075, -.36484967827935958644932861307814066351067871619058) {7};
\node at (.21093242412093931363351251934884378021978648378953, 1.2764824302262961524609825137197259270450311761020, -.86484967827935958644932861307814066351067871619060) {8};
\node at (-.26811812545734779130689058944601519996404239063795, .42255219728220816465295019618437087043498939248147, -1.0681074658615250936386579889886538538269692085900) {9};
\end{tikzpicture}

\begin{tikzpicture}
\tikzset{
    table/.style={
        matrix of nodes,
        row sep=-\pgflinewidth,
        column sep=-\pgflinewidth,
        nodes={
            rectangle,
            draw=black,
            align=center
        },
        minimum height=1.5em,
        text depth=1ex,
        text height=1.5ex,
        nodes in empty cells,
        every even row/.style={
            nodes={fill=gray!20}
        },
        column 1/.style={
            nodes={text width=1em,font=\bfseries}
        },
        row 1/.style={
            nodes={
                fill=black,
                text=white,
                font=\bfseries
            }
        }
    }
}
 \pgfgettransform\mytrafo
\matrix (first) [table,text width=10em,execute at begin cell=\pgfsettransform\mytrafo]
{
& X & Y & Z\\
1 & $0$ & $0$ & $0$ \\
2 & $1$ & $0$ & $0$ \\
3 & $1/2$ & $\sqrt{3}/2$ & $0$ \\
4 & $1/2$ & {\tiny $-(9\alpha^{3}+3\alpha^{2}-9\alpha-2)/2\sqrt{3}$} & {\tiny $\sqrt{3-3\alpha^{2}}(9\alpha^{3}+3\alpha^{2}-3\alpha-1)/4$} \\
5 & $(27\alpha^{3}-24\alpha-2)/3$ & $-\sqrt{3}(9\alpha^{3}-8\alpha-3)/7$ & {\tiny $\sqrt{3-3\alpha^{2}}(9\alpha^{3}+3\alpha^{2}-3\alpha-1)/14$} \\
6 & $(3-3\alpha)/4$ & $-\sqrt{3}(-1+\alpha)/4$ & $-\sqrt{3-3\alpha^{2}}/2$ \\
7 & $9\alpha^{2}/8-3\alpha/4+1/8$ & $\sqrt{3}(3\alpha^{2}+2\alpha-1)/8$ & $\sqrt{3-3\alpha^{2}}(3\alpha-1)/4$ \\
8 & $(1-3\alpha)/4$ & $-\sqrt{3}(\alpha-3)/4$ & $-\sqrt{3-3\alpha^2}/2$ \\
9 & $(-27\alpha^{3}+24-5)/14$ & {\tiny $-\sqrt{3}(9\alpha^{3}-8\alpha-3)/14$ }& {\tiny $\sqrt{3-3\alpha^{2}}(81\alpha^{3}-51\alpha-6)/14$} \\
};
\end{tikzpicture}
\begin{figure}[!htb]
\centering
\caption{Soifer Graph:  $\alpha=0.05209\ldots$  is the smallest root of the equation $27z^{4}+18z^{3}-18z+1$}
\label{}
\end{figure}

Select the vertices $A_1$, $A_2$ and $A_3$ as the vertices of an unit equilateral triangle as shown.
Then pick $A_4$ such that $\|A_4-A_1\|=\|A_4-A_2\|=1$, one degree of freedom left. Similarly, let
$A_5$ such that $\|A_5-A_1\|=\|A_5-A_3\|=1$, pick $A_6$ such that $\|A_6-A_2\|=\|A_6-A_3\|=1$, let $A_7$ such that $\|A_7-A_2\|=\|A_7-A_6\|=1$, let $A_8$ such that $\|A_8-A_3\|=\|A_8-A_6\|=1$, and finally $A_9$ such that $\|A_9-A_6\|=\|A_9-A_7\|=1$. At this point we have six parameters and there are five unit edges unaccounted for:
$\|A_4-A_5\|$, $\|A_4-A_7\|$, $\|A_4-A_9\|$, $\|A_5-A_8\|$, and $\|A_8-A_9\|$. The system is underdetermined;
the solution presented above is the simplest we could find.
\end{section}

\section{\bf Conclusions and directions of future research}

The problem of computing the Euclidean dimension of a given graph is an interesting and difficult one.
In this paper, we determined this quantity for twelve well known graphs. Admittedly, we focused on rather
small graphs, but even so, some of the embeddings were rather challenging to find.

While a solution to the general problem seems hopeless at the present time, it is reasonable to expect that
one can still uncover interesting results. One may try to find a better upper bound
than the estimate of Maehara and R\"{o}dl mentioned in Theorem \ref{MR}. One such idea is sketched below.

The {\it chromatic number} of a graph $G$, denoted $\chi(G)$, is the minimum number of colors that
can be assigned to the vertices of $G$ such that no two adjacent vertices receive the same color.
Recall than $\Delta(G)$ denotes the maximum degree over all vertices of $G$.

Brooks \cite{Brooks} proved that $\chi(G)\le \Delta(G)$ unless $G$ is a complete graph or an odd cycle, case in which $\chi(G)=\Delta+1$. We propose the following
\begin{conj}
For every graph $G$ we have that $dim(G)\le 2\chi(G)$.
\end{conj}

As per Brooks' theorem, this would represent an improvement of Maehara and R\"{o}dl's bound.
The connection between the chromatic number and the Euclidean dimension is also interesting for another reason:
the Hadwiger-Nelson problem.

In the early 1950s Hadwiger and Nelson asked what is the minimum number of colors needed to color the plane
such that no two points unit distance apart are colored identically. It is known that this number is $4,\, 5,\, 6$ or $7$.

By a result of de Brujin and Erd\H{o}s \cite {dBE}, this is equivalent to asking what is the maximum chromatic
number of a graph whose Euclidean dimension is $2$. There has been no progress on this problem for more than 65 years.
\newpage
\thispagestyle{empty}

\end{document}